\documentclass[reqno]{amsart}

\usepackage{hyperref}

\newtheorem{theorem}{Theorem}
\newtheorem{lemma}[theorem]{Lemma}
\newtheorem{proposition}[theorem]{Proposition}

\newtheorem{example}[theorem]{Example}
\theoremstyle{exercise}

\theoremstyle{definition}
\newtheorem{definition}[theorem]{Definition}
\usepackage{graphicx}
\usepackage{pstricks, enumerate, pst-node, pst-text, pst-plot}
\theoremstyle{remark}
\newtheorem{remark}[theorem]{Remark}

\numberwithin{equation}{section}

\newcommand{\intav}[1]{\mathchoice {\mathop{\vrule width 6pt height 3 pt depth  -2.5pt
\kern -8pt \intop}\nolimits_{\kern -6pt#1}} {\mathop{\vrule width
5pt height 3  pt depth -2.6pt \kern -6pt \intop}\nolimits_{#1}}
{\mathop{\vrule width 5pt height 3 pt depth -2.6pt \kern -6pt
\intop}\nolimits_{#1}} {\mathop{\vrule width 5pt height 3 pt depth
-2.6pt \kern -6pt \intop}\nolimits_{#1}}}

\newcommand{\intavl}[1]{\mathchoice {\mathop{\vrule width 6pt height 3 pt depth  -2.5pt
\kern -8pt \intop}\limits_{\kern -6pt#1}} {\mathop{\vrule width 5pt
height 3  pt depth -2.6pt \kern -6pt \intop}\nolimits_{#1}}
{\mathop{\vrule width 5pt height 3 pt depth -2.6pt \kern -6pt
\intop}\nolimits_{#1}} {\mathop{\vrule width 5pt height 3 pt depth
-2.6pt \kern -6pt \intop}\nolimits_{#1}}}

\begin{document}

\title[Non-uniformly hyperbolic horseshoes]{Fractal geometry of non-uniformly hyperbolic horseshoes}

\author[Carlos Matheus]{Carlos Matheus}
\address{Universit\'e Paris 13, Sorbonne Paris Cit\'e, LAGA, CNRS (UMR 7539), F-93430, Villetaneuse, France.}
\email{matheus@impa.br}
\thanks{I'm thankful to Idris Assani for the kind invitation to write this article. Also, I'm grateful to my coauthors J. Palis and J.-C. Yoccoz for the pleasure of working with them on the project at the origin of this article, to the anonymous referees for their useful comments, and to Coll\`ege de France for the hospitality during the preparation of this text. Last but not least, the author was partially supported by the Balzan research project of J. Palis.}


\date{\today}


\begin{abstract}
This article is an expanded version of some notes for my talk at the ``Ergodic Theory and Dynamical Systems Workshop'' (from March 22 to March 25, 2012) held at the Department of Mathematics of the University of North Carolina at Chapel Hill. In the aforementioned talk, it was discussed some recent results on the fractal geometry of certain objects -- \emph{non-uniformly hyperbolic horseshoes} -- constructed by Jacob Palis and Jean-Christophe Yoccoz in their recent \emph{tour-de-force} work around \emph{heteroclinic bifurcations of surface diffeomorphisms}.

The goal of the present article is two-fold. The first part will be a modest survey on the history of the study of homoclinic/heteroclinic bifurcations of surface diffeomorphisms: starting from the seminal works of Henri Poincar\'e on
Celestial Mechanics, we will recall some landmark results on bifurcations of homoclinic/heteroclinic tangencies
associated to \emph{uniformly hyperbolic horseshoes}. Our discussions in this (first) part will culminate with a brief presentation of the scheme of the
proof of the main theorems of J. Palis and J.-C. Yoccoz on this subject. In particular, we will highlight the main features of the so-called non-uniformly hyperbolic horseshoes, an object at the heart of the work of J. Palis and J.-C. Yoccoz. Then, the second part will be a sort of research announcement where we will discuss some results obtained in collaboration with J. Palis and J.-C. Yoccoz on the (fractal) geometry of non-uniformly hyperbolic horseshoes.
\end{abstract}

\maketitle

\tableofcontents


\section{Part I -- a survey on homoclinic/heteroclinic bifurcations}

In his seminal work (in 1890) on Celestial Mechanics, Henri Poincar\'e \cite{P} emphasized the relevance of the concept of \emph{homoclinic} orbits in Dynamical Systems by stating:

\begin{center}
``... \emph{rien n'est plus propre \`a nous donner une id\'ee de la complication de tous les probl\`emes de dynamique} ...''
\end{center}
\begin{center}
(in a free translation to English: ``... \emph{nothing is more adequate to give us an idea of the complexity of all problems in dynamics} ...'')
\end{center}

In fact, the history behind the introduction of this notion is fascinating: in a few words, H. Poincar\'e submitted a first version \cite{P-Acta} of his work to a competition in honor of G. Mittag-Leffler and financially support by the king Oscar of Sweden, but, after some comments of L. Phragm\'en, it was discovered a mistake in part of his text related to the presence of homoclinic orbits. For nice accounts in English and French (resp.) on this beautiful chapter of the history of Dynamical Systems, see \cite{B-G} and \cite{Y} (resp.).

In modern language, we define an homoclinic orbit as follows. Given a diffeomorphism $f:M\to M$ of a compact (boundaryless) manifold $M$, denote by $f^n =
\underbrace{f \circ\dots\circ f}_{n}$ the $n$-th iterate of $f$, $n\in\mathbb{Z}$. Let $p\in M$ be a periodic point of $f$ with minimal period $k$, i.e., $f^{k}(p)=p$ and $k\in\mathbb{N}$ is minimal with this property. We say that the orbit $\{f^n(q): n\in\mathbb{Z}\}$ of a point $q\neq p$ is \emph{homoclinic} to $p$ whenever $f^{nk}(q)\to p$ as $n\to\pm\infty$, that is, the orbit of $q$ is accumulates the orbit of the periodic point $p$ both in the past and the future.

Similarly, given two \emph{distinct}\footnote{Here, we mean that $p_1$ and $p_2$ belong to distinct \emph{orbits}, i.e., $p_2\neq f^n(p_1)$ for all $n\in\mathbb{N}$.} periodic points $p_1, p_2\in M$ with (minimal) periods $k_1, k_2$ (resp.), we say that the orbit of a point $q\neq p_1, p_2$ is \emph{heteroclinic} to $p_1$ and $p_2$ whenever $f^{n k_1}(q)\to p_1$ as $n\to-\infty$ and
$f^{n k_2}(q)\to p_2$ as $n\to+\infty$.

George Birkhoff was one of the first to confirm the predictions of H. Poincar\'e on homoclinic orbits by proving in 1935 that, in general, one can find periodic orbits of very high period near homoclinic orbits.

Later on, by taking as a source of inspiration the works of G. Birkhoff on homoclinic orbits, and Cartwright and Littlewood \cite{CL}, \cite{L1} and \cite{L2}, and Levinson \cite{Le} on differential equations similar to the Van der Pol equation\footnote{A differential equation steaming from Engineering problems related to nonlinear oscillators and radio waves.}, Steve Smale proposed in 1967 a geometrical model currently referred to as \emph{Smale's horseshoe} explaining in a very satisfactory way the mechanism responsible for the dynamical complexity near a general homoclinic orbit.

In the subsection below, we will quickly revisit some features of Smale's horseshoe as a paradigm of \emph{hyperbolic set} of a dynamical system. The basic references for historical and mathematical details on the next three subsections is the classical book \cite{PT} of J. Palis and F. Takens.


\subsection{Transverse homoclinic orbits and Smale's horseshoes} Let $f:M\to M$ be a $C^k$ diffeomorphism, $k\geq 1$ and let $p\in M$ be a periodic point of $f$. For sake of simplicity\footnote{This can be achieved by replacing $f$ by some iterate $f^k$, and, as far as the discussion of this subsection is concerned, this replacement has no serious effect.}, let's assume that $p$ is a fixed point, i.e., $f(p)=p$. The \emph{stable} and \emph{unstable} sets of $p$ are
$$W^s(p):=\{q\in M: f^n(q)\to p \textrm{ as } n\to+\infty\} \quad \textrm{ and } \quad W^u(p):=\{q\in M: f^n(q)\to p \textrm{ as } n\to-\infty\}$$

In this notation, $q$ is homoclinic to $p$ if and only if $q\in (W^s(p)\cap W^u(p))-\{p\}$, and $q$ is heteroclinic to $p_1$ and $p_2$ if and only if
$q\in (W^u(p_1)-\{p_1\})\cap (W^s(p_2)-\{p_2\})$.

For a generic $f$, the fixed point is \emph{hyperbolic}, i.e., the differential $df(p):T_pM\to T_pM$ is a linear map without eigenvalues of norm $1$. In this case, denote by $E^s$, $E^u$, the stable and unstable subspaces of $df(p)$, i.e., the generalized eigenspaces of $df(p)$ associated to the eigenvalues of norm strictly smaller, resp. larger, than $1$. Then, by the \emph{stable manifold theorem}\footnote{Cf. Appendix 1 of \cite{PT}.}, the stable and unstable sets of $p$ (i.e., $W^s(p)$ and $W^u(p)$) are injectively \emph{immersed} $C^k$ submanifolds of $M$ of dimension $s$, $u$, where $s=\textrm{dim}(E^s)$, $u=\textrm{dim}(E^u)$.

We say that $q$ is a \emph{transverse homoclinic orbit} to a hyperbolic fixed point $p$ when the stable and unstable manifolds of $p$ intersect transversally at $q\neq p$, that is,  $q\in (W^s(p)\cap W^u(p))-\{p\}$ and $T_qM=T_q W^s(p)\oplus T_qW^u(p)$. By transversality theory (or more precisely, Kupka-Smale's theorem), for a generic $f$, all homoclinic orbits to periodic points are transverse.

The fundamental picture discovered by S. Smale near a transverse homoclinic orbits to hyperbolic fixed points is the following.

\begin{figure}[htb!]
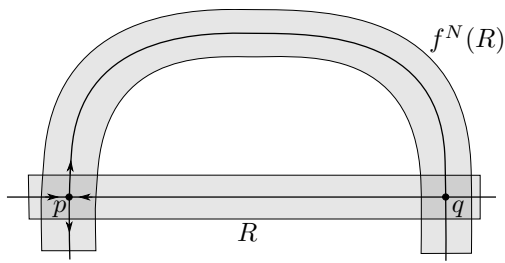
\caption{Smale's horseshoe.}\label{f.1}
\end{figure}


In a nutshell, this picture means that, near a point $q$ which is transverse homoclinic to a hyperbolic fixed point $p$, one can find a rectangle $R$ containing $p$ and $q$ such that some iterate $F=f^N$ of $f$ maps $R$ in the ``horseshoe''-shaped region $f^N(R)$ shown above. Moreover, the picture was drawn to convince the reader that the action of the differential $dF$ of $F$ on $R$ uniformly contracts any almost horizontal direction and uniformly expand any almost vertical direction.

Using these facts, S. Smale proved that the \emph{maximal invariant set} $\Lambda:=\bigcap\limits_{n\in\mathbb{Z}}f^{nN}(R) = \bigcap\limits_{n\in\mathbb{Z}}F^{n}(R)$ consisting of all points in $R$ whose orbit under $F$ never escapes $R$ is a \emph{hyperbolic set}, that is, there are constants $C>0$, $0<\lambda<1$ and a splitting $T_xM=E^s(x)\oplus E^u(x)$ for each $x\in\Lambda$ such that:
\begin{itemize}
\item the splitting is $dF$-invariant: $dF(x)\cdot E^s(x)=E^s(F(x))$ and $dF(E^u(x))=E^u(F(x))$;
\item $E^s$ is uniformly contracted and $E^u$ is uniformly expanded: $\|dF^n(x)\cdot v^s\|, \|dF^{-n}(x)\cdot v^u\|\leq C\lambda^n\|v\|$ for all $n\geq 0$, $v^s\in E^s(x)$,
$v^u\in E^u(x)$, where $\|.\|$ is a norm associated to some choice of Riemannian metric on $M$.
\end{itemize}

\begin{remark}
In the case of our picture above, there is no mystery behind the choice of the splitting: $E^s(x)$ is an almost horizontal direction and $E^u(x)$ is an almost vertical direction.
\end{remark}

Furthermore, by using the hyperbolicity of the set $\Lambda$, S. Smale showed that the dynamics of $F$ restricted to $\Lambda$ is \emph{topologically conjugated to Bernoulli shift in two symbols}, that is, there exists a homeomorphism $h:\Lambda\to\Sigma:=\{0,1\}^{\mathbb{Z}}$ such that $h(F(x))=\sigma(h(x))$ where $\sigma:\Sigma\to\Sigma$ is given by $\sigma((a_i)_{i\in\mathbb{Z}})=(a_{i+1})_{i\in\mathbb{Z}}$. In other words, the dynamics of $F|_{\Lambda}$ can be modeled by a Markov process.

Among the several striking consequences of S. Smale's results, we observe that the set of periodic orbits of $F$ is dense in $\Lambda$ and the dynamical system
$F|_{\Lambda}$ is \emph{sensitive to initial conditions}\footnote{That is, two nearby distinct points tend to get far apart after an appropriate number of iterations of the dynamics.} simply because the same is true for the Bernoulli shift $\sigma$! In particular, S. Smale's results allow one to recover the result of G. Birkhoff (mentioned above) that the transverse homoclinic point $q$ of the hyperbolic periodic point $p$ is accumulated by periodic orbits of $f$ of arbitrarily high period.

By obvious reason, the maximal invariant set $\Lambda$ was baptized \emph{horseshoe} by S. Smale. Partly motivated by this, we introduce the following concepts:

\begin{definition}\label{d.hyp-set}We say that a compact subset $\Lambda\subset M$ is a \emph{hyperbolic set} of a diffeomorphism $f:M\to M$ if
\begin{itemize}
\item $\Lambda$ is $f$-invariant, that is, $f(\Lambda)=\Lambda$;
\item there are constants $C>0$, $0<\lambda<1$ and a splitting $T_xM=E^s(x)\oplus E^u(x)$ for each $x\in\Lambda$ with:
\begin{itemize}
\item $df(x)\cdot E^s(x)=E^s(f(x))$ and $df(E^u(x))=E^u(f(x))$;
\item $\|df^n(x)\cdot v^s\|, \|df^{-n}(x)\cdot v^u\|\leq C\lambda^n\|v\|$ for all $n\geq 0$, $v^s\in E^s(x)$,
$v^u\in E^u(x)$, where $\|.\|$ is a norm associated to some choice of Riemannian metric on $M$.
\end{itemize}
\end{itemize}
\end{definition}

In other words, $\Lambda\subset M$ is a hyperbolic set of a diffeomorphism $f:M\to M$ whenever $\Lambda$ is $f$-invariant, and the differential $df$ completely decomposes $T_{\Lambda}M$ into two $df$-equivariant subbundles $E^s$ and $E^u$ such that $E^s$ is a \emph{stable} subbundle (that is, it is forwardly contracted by $df$) and $E^u$ is an \emph{unstable} subbundle (that is, it is backwardly contracted by $df$).

\begin{example} The orbit $\Lambda=\{p,\dots, f^{k-1}(p)\}$ of a hyperbolic period point $p$ of period $k$ is a trivial (i.e., finite) hyperbolic set, while Smale's horseshoes are non-trivial (i.e., infinite) hyperbolic sets.
\end{example}

One of the key features of hyperbolic sets is the fact that the infinitesimal information on the structure of $df$ over a hyperbolic set $\Lambda$ imposes a certain number of global geometrical constraints on the dynamics of $f$ on $\Lambda$. For example, given $x\in M$, denote by $W^s(x)=\{y\in M: \textrm{dist}(f^n(y), f^n(x))\to 0 \textrm{ as }n\to+\infty\}$ and $W^u(x)=\{y\in M: \textrm{dist}(f^n(y), f^n(x))\to 0 \textrm{ as }n\to-\infty\}$ the stable and unstable sets of $x$. In general, the stable and unstable sets of an arbitrary point of an arbitrary diffeomorphism may have a \emph{wild} geometry, such as fractal sets. On the other hand, as we already mentioned, it is known that the stable and unstable sets of hyperbolic periodic points are injectively immersed submanifolds thanks to the stable manifold theorem. In other words, the geometry of stable and unstable sets improves under appropriate hyperbolicity conditions, and, as it turns out, it is possible to generalize the stable manifold theorem to show that the stable and unstable sets of any point in a hyperbolic set has well-behaved stable and unstable sets:

\begin{theorem}[Generalized stable manifold theorem] Let $\Lambda\subset M$ be a hyperbolic set of a $C^k$-diffeomorphism $f:M\to M$, $k\geq 1$. Then, the stable set $W^s(x)$ of any $x\in\Lambda$ is an injectively immersed $C^k$-submanifold of dimension $\textrm{dim}(E^s(x))$ and, for all sufficiently small $\varepsilon>0$, the local stable set 
$$W^s_{loc}(x)=\{y\in W^s(x): \textrm{dist}(f^n(y),f^n(x))\leq\varepsilon\}$$ 
is a $C^k$ embedded disk in $W^s(x)$ of dimension $\textrm{dim}(E^s(x))$.
Also, the stable sets $W^s(x)$ depend continuously on $x\in\Lambda$ and $f$. Furthermore, the map $\Lambda\ni x\mapsto W^s_{loc}(x)\subset M$ is continuous. 
\end{theorem}

Another way of phrasing the previous theorem is: given a hyperbolic set $\Lambda$, the family of stable sets $W^s(x)$ of points $x\in\Lambda$ form a continuous \emph{lamination} with $C^k$ leaves.

Actually, this is not the full statement of the generalized stable manifold theorem. For more complete statements see Appendix 1 of \cite{PT} and references therein (especially \cite{HPS}). 

Coming back to the discussion of Smale's horseshoes, it turns out that they are not arbitrary hyperbolic sets in the sense that they fit the following definitions:

\begin{definition}\label{d.basic-set} A set $\Lambda\subset M$ is a \emph{basic set} of a diffeomorphism $f:M\to M$ if $\Lambda$ is an \emph{infinite} hyperbolic set such that
\begin{itemize}
\item $\Lambda$ is \emph{transitive}, i.e., there exists $x\in\Lambda$ whose orbit $\{f^n(x)\}_{n\in\mathbb{Z}}$ is dense in $\Lambda$;
\item $\Lambda$ is \emph{locally maximal}, i.e., there exists a neighborhood $U$ of $\Lambda$ such that the maximal invariant set $\bigcap\limits_{n\in\mathbb{Z}}f^n(U)$ of $U$ coincides with $\Lambda$, that is, $\bigcap\limits_{n\in\mathbb{Z}}f^n(U)=\Lambda$.
\end{itemize}
\end{definition}

The notion of basic set is natural in our setting because the transitivity and local maximality properties allow one to show that the hyperbolicity of $\Lambda$ is a robust property in the sense that the set $\Lambda_g:=\bigcap\limits_{n\in\mathbb{Z}}g^n(U)$ (called \emph{continuation} of $\Lambda$) is a hyperbolic set whenever $g$ is sufficiently $C^k$-close to $f$. See, e.g., the book \cite{PT} and the references therein for more details.

\begin{definition}\label{d.UH-horseshoe}A set $\Lambda\subset M$ is a \emph{(uniformly hyperbolic) horseshoe} of a diffeomorphism $f:M\to M$ if $\Lambda$ is a totally disconnected basic set of $f$ of \emph{saddle-type}, i.e.,
both subbundles $E^s$ and $E^u$ appearing in Definition \ref{d.hyp-set} are non-trivial.
\end{definition}

Concerning this definition, let us mention that in these notes we'll focus exclusively on saddle-type hyperbolic sets because they are the most relevant for the study of homoclinic/heteroclinic bifurcations (as in this context we need, by definition, both stable and unstable manifolds). However, it is worth it to point out that the dynamics of attractors and/or repellors (i.e., the situations where either $E^s$ or $E^u$ is trivial) is also very exciting and it is not surprising that they have a vast literature dedicated to them (see for instance the book \cite{PT} and references therein).

From the qualitative point of view, a uniformly hyperbolic horseshoe $\Lambda$ of a diffeomorphism $f$ behaves exactly as a Smale's horseshoe near a transverse homoclinic orbit. For instance, it is possible to show that the restriction $f$ to $\Lambda$ is topologically conjugated to a Markov shift of finite type. In particular, $\Lambda$ is topologically a \emph{Cantor set}\footnote{A non-empty compact totally disconnected set that is perfect (i.e., without isolated points).}, and, despite the fact that the dynamics of $f|_{\Lambda}$ is \emph{chaotic} (e.g., in the sense that it is sensitive to initial conditions), one can reasonably understand $f|_{\Lambda}$ because it topologically modeled by a \emph{Markov process}.\footnote{Actually, it is possible to prove that $f|_{\Lambda}$ also has plenty of interesting properties from the statistical (\emph{ergodic}) point of view: for example, $\Lambda$ supports several ergodic $f$-invariant probabilities coming from the so-called \emph{thermodynamical formalism} of R. Bowen, D. Ruelle and Y. Sinai. See e.g. \cite{B} for more nice account on this subject.}

Therefore, we can declare that the local dynamics near transverse homoclinic orbits, or, more generally, uniformly hyperbolic horseshoes, is well-understood, and hence we can start the discussion of the local dynamics near homoclinic \emph{tangencies} (i.e., non-transverse homoclinic orbits).


\subsection{Homoclinic tangencies and Newhouse phenomena}\label{ss.bifurcations-Newhouse} Let $K$ be a (uniformly hyperbolic) horseshoe of a $C^k$, $k\geq 2$, diffeomorphism $f:M\to M$ of a compact \emph{surface} (two-dimensional manifold) $M$ possessing a periodic point $p\in K$ associated to a \emph{quadratic homoclinic tangency}, that is, the stable and unstable manifolds (curves in our current setting) $W^s(p)$ and $W^u(p)$ of $p$ meet tangentially at a point $q\in (W^s(p)\cap W^u(p))-K$ and the \emph{order of contact} between $W^s(p)$ and $W^u(p)$ at $q$ is $1$, that is, the curves $W^s(p)$ and $W^u(p)$ are tangent at $q$ but their \emph{curvatures} differ at $q$.

The main geometrical features of a quadratic homoclinic tangency are captured by the picture in Figure \ref{f.2} below.

\begin{figure}[htb!]
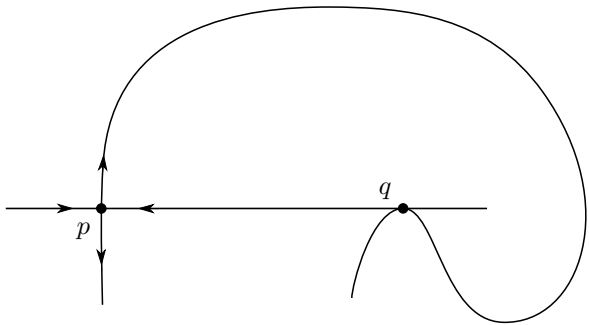
\caption{Quadratic homoclinic tangency associated to a periodic point in a horseshoe.}\label{f.2}
\end{figure}

For sake of simplicity, we'll assume that there are two neighborhoods $U$ of the horseshoe $K$ and $V$ of the homoclinic orbit $\mathcal{O}(q):=\{f^n(q):n\in\mathbb{Z}\}$ of $q$ such that
$$\bigcap\limits_{n\in\mathbb{Z}} f^n(U\cup V)=K\cup \mathcal{O}(q).$$
In other words, we'll suppose that the local dynamics of $f$ on $U\cup V$ consist precisely of the horseshoe $K$ and the homoclinic orbit of tangency $\mathcal{O}(q)$, that is, locally (on $U\cup V$) the interesting dynamical phenomena come exclusively from the horseshoe and $\mathcal{O}(q)$. See Figure~\ref{f.3} below.

\begin{figure}[hbt!]
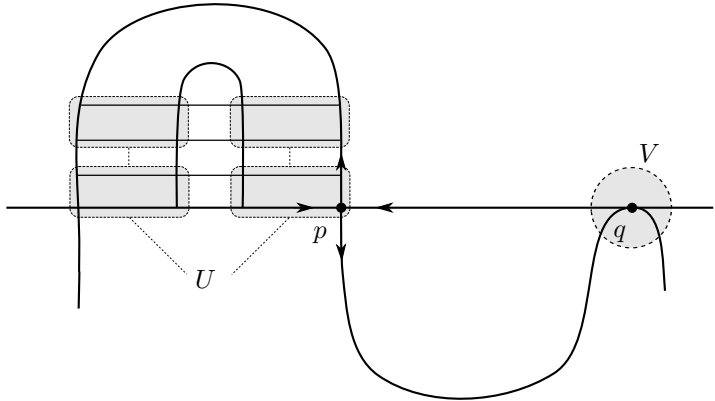
\caption{Localization of the dynamics via the neighborhoods $U$ (of the horseshoe) and $V$ (of the tangency).}\label{f.3}
\end{figure}

Note that the maximal invariant set
$$\Lambda_f:=\bigcap\limits_{n\in\mathbb{Z}} f^n(U\cup V)$$
capturing the local dynamics of $f$ on $U\cup V$ is not a hyperbolic set. Indeed, it is not hard to convince oneself that the natural candidate for the stable $E^s(q)$, resp. unstable $E^u(q)$, direction at $q$ in Definition \ref{d.hyp-set} is the $1$-dimensional direction $T_qW^s(p)$, resp. $T_qW^u(p)$. However, since $W^s(p)$ and $W^u(p)$ meet tangentially at $q$, one would have $E^s(q)=T_qW^s(p)=T_qW^u(p)=E^u(q)$, so that the condition $T_qM = E^s(q)\oplus E^u(q)$ in Definition \ref{d.hyp-set} is never fulfilled.

On the other hand, since $\Lambda_f=K\cup\mathcal{O}(q)$ and the single orbit $\mathcal{O}(q)$ is solely responsible for the non-hyperbolicity of $\Lambda_f$, we still completely understand the local dynamics of $f$ on $U\cup V$.

Now, let's try to understand the local dynamics on $U\cup V$ of a $C^k$-diffeomorphism $g:M\to M$ $C^k$-\emph{close} to $f$. Consider
$\mathcal{U}$ a sufficiently small $C^k$-neighborhood of $f$ such that the dynamically relevant objects in Figure \ref{f.3} above admit a \emph{continuation} for any
$g\in\mathcal{U}$: more precisely, we select $\mathcal{U}$ so that, for any $g\in\mathcal{U}$, the maximal invariant set
$$K_g=\bigcap\limits_{n\in\mathbb{Z}}g^{n}(U)$$
is a (uniformly hyperbolic) horseshoe (cf. the paragraph after Definition~\ref{d.basic-set}), the periodic point $p$ has a continuation into a nearby (hyperbolic) periodic point $p_g$ of $g$, and the compact curve
$c^s(f)$, resp. $c^u(f)$, inside the stable, resp. unstable, manifold $W^s(p)$, resp. $W^u(p)$ containing $p$ and $q$ and crossing $V$ has a continuation into a nearby compact curve $c^s(g)$, resp. $c^u(g)$, in the stable, resp. unstable, manifold of
$p_g$ crossing $V$.

Using these dynamical objects associated to $g\in\mathcal{U}$, we can organize the parameter space $\mathcal{U}$ by writing $\mathcal{U}=\mathcal{U}_-\cup \mathcal{U}_0\cup \mathcal{U}_+$ where
\begin{itemize}
\item $g\in\mathcal{U}_-$ whenever $c^s(g)$ and $c^u(g)$ don't intersect;
\item $g\in\mathcal{U}_0$ whenever $c^s(g)$ and $c^u(g)$ have a quadratic tangency at a point $q_g$ in $V$;
\item $g\in\mathcal{U}_+$  whenever $c^s(g)$ and $c^u(g)$ have two transverse intersection points in $V$.
\end{itemize}

Since $q$ corresponds to a quadratic tangency of $f$, we have that $\mathcal{U}_0$ is a codimension $1$ hypersurface dividing $\mathcal{U}$ into the two connected open subsets $\mathcal{U_-}$ and $\mathcal{U}_+$. The picture below illustrates the decomposition $\mathcal{U}=\mathcal{U}_-\cup \mathcal{U}_0\cup \mathcal{U}_+$ of the parameter space and the features on phase space of the elements of $\mathcal{U}_-$, $\mathcal{U}_0$ and $\mathcal{U}_+$.

\begin{figure}[hbt!]
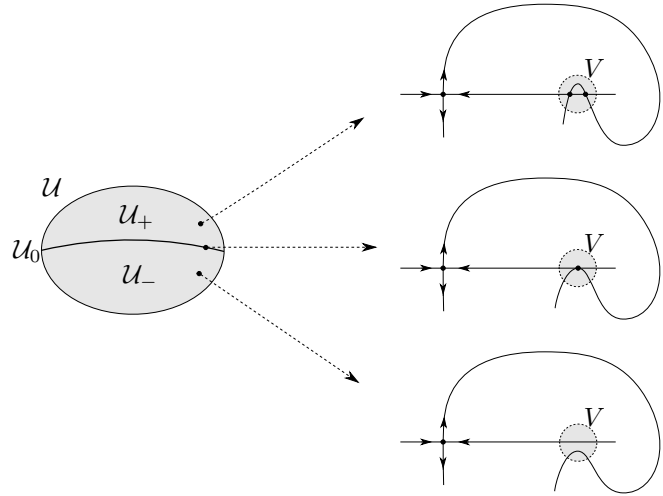
\caption{Organization of the parameter space $\mathcal{U}$.}\label{f.4}
\end{figure}

From the (local) dynamical point of view, the regions $\mathcal{U}_-$ and $\mathcal{U}_0$ of the parameter space $\mathcal{U}$ are not particularly interesting: in fact, by inspecting the definitions, it is not hard to show that
\begin{itemize}
\item $\Lambda_g=K_g$ for any $g\in\mathcal{U}_-$, and 
\item $\Lambda_g=K_g\cup \mathcal{O}(q_g)$ for any $g\in\mathcal{U}_0$.
\end{itemize}

In other words, all potentially new dynamical phenomena come from $\mathcal{U}_+$, that is, after non-trivially unfolding the quadratic tangency associated to diffeomorphisms in $\mathcal{U}_0$. 

Here, it may be tempting to try to understand the local dynamics of \emph{all} $g\in\mathcal{U}_+$. However, after the seminal works of Sheldon Newhouse \cite{N1}, \cite{N2}, \cite{N3}, one knows that it is not reasonable to try to control $\Lambda_g$ for \emph{all} $g\in\mathcal{U}_+$ because of a mechanism nowadays called \emph{Newhouse phenomena}. 

More precisely, it is clear from Figure \ref{f.4} above that by unfolding a quadratic tangency to get an element $g\in\mathcal{U}_+$, we end up with 
some horseshoes near the region $V$: indeed, the presence of transverse homoclinic orbits in $V$ of $g\in\mathcal{U}_+$ implies the existence of some horseshoes by the discussion of the previous subsection. In particular, this naive argument seems to ensure that $\Lambda_g$ is always hyperbolic. However, S. Newhouse noticed that by unfolding the quadratic tangency associated to $p$ we may create other tangencies nearby, that is, we may ``accidentally'' lose the hyperbolicity just created in view of the incompatibility between tangencies and hyperbolicity. 

In fact, as we're going to explain in a moment, this ``accidental'' formation of new tangencies happens especially when the horseshoe $K$ containing $p$ is thick (fat). In this direction, the first step is to reduce the detection of tangencies for diffeomorphisms of $2$-dimensional manifolds to the ($1$-dimensional) problem of understanding the intersection of two Cantor sets in $\mathbb{R}$.

\subsubsection{Persistence of tangencies and intersections of Cantor sets}\label{sss.tangencies-Cantors}

Starting from $f\in\mathcal{U}_0$, we consider an extension\footnote{Such an extension exists by the results of W. de Melo \cite{dM} and it heavily depends on the fact that $f$ is a diffeomorphism of a $2$-\emph{dimensional} manifold.} $\mathcal{F}^s(f)$, resp., $\mathcal{F}^u(f)$, of the stable, resp. unstable, laminations of $K$ to the neighborhood $U$ of $K$. From the fact that $q$ is a homoclinic quadratic tangency associated to $p$, one deduces that the foliations $\mathcal{F}^s$ and $\mathcal{F}^u$ meet tangencially along a curve $\ell=\ell(f)$ called the \emph{line of tangencies}. Using 
$\ell(f)$ as an auxiliary curve, we can study the formation of tangencies for $g\in\mathcal{U}_+$ (i.e., after unfolding the quadratic tangency of $f$) as follows. One considers the (local) Cantor sets $W_{loc}^s(p)\cap K$ and $W_{loc}^u(p)\cap K$, and, by using the holonomy of the stable, resp. unstable, foliations $\mathcal{F}^s(f)$, resp. $\mathcal{F}^u(f)$, one can diffeomorphically map $W_{loc}^s(p)\cap K$ and $W_{loc}^u(p)\cap K$ into Cantor sets $K^s\subset\ell$ and $K^u\subset \ell$ by sending $x\in W_{loc}^s(p)\cap K$, resp. $x\in W_{loc}^s(p)\cap K$, to $\pi_f^s(x)=y\in W^s(x)\cap\ell$, resp. $\pi_f^u(x)=y\in W^u(p)\cap\ell$. Note that, by definition, the intersection of $K^s$ and $K^u$ corresponds to all tangencies between the stable and unstable laminations of the horseshoe $K$ near $V$, that is, by our assumptions, $K^s\cap K^u=\{q\}$. Pictorially, the discussion of this paragraph can be illustrated by the following figure:

\begin{figure}[hbt!]
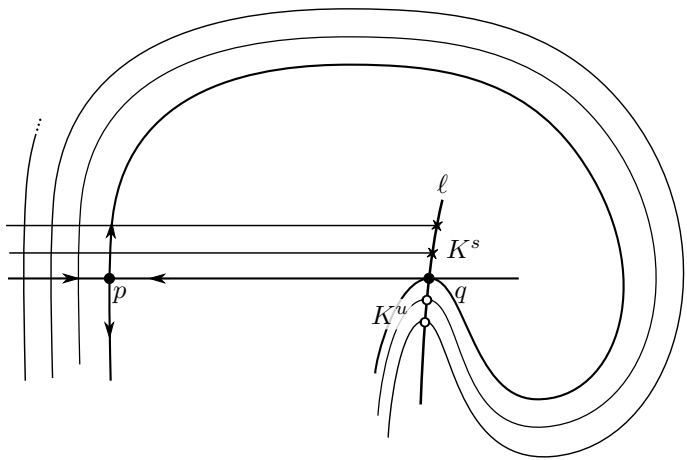
\caption{The line of tangencies $\ell$ and the Cantor sets $K^s$ and $K^u$ for $f\in\mathcal{U}_0$: the crosses are points in $K^s$ and the dots are points in $K^u$.}\label{f.5}
\end{figure}

Now, let us perturb this picture by unfolding the tangency to get some $g\in\mathcal{U}_+$. It is possible to show that this picture admits a natural continuation because the relevant dynamical objects have continuations\footnote{We'll comment more on this in Subsection \ref{ss.gap-lemma}, but for now let's postpone this ``continuity'' discussion.}, that is, one can extend the stable and unstable laminations of $K_g$ to stable and unstable foliations 
$\mathcal{F}^s(g)$ and $\mathcal{F}^u(g)$ of a neighborhood $U$ (close to $\mathcal{F}^s(f)$ and $\mathcal{F}^u(f)$), and one can use them to define a line of tangencies 
$\ell(g)$ (``close'' to $\ell(f)$) containing two Cantor sets $K^s(g)$ and $K^u(g)$ (close to $K^s$ and $K^u$) defined as the images of the Cantor sets $W_{loc}^s(p_g)\cap K_g$ and $W_{loc}^u(p_g)\cap K_g$ via the stable and unstable holonomies $\pi^s_g$ and $\pi^u_g$. Also, the intersection $K^s(g)\cap K^u(g)$ of the Cantor sets $K^s(g)$ and $K^u(g)$ still  accounts for all tangencies between the stable and unstable laminations of the horseshoe $K_g$. In summary, we get the following local dynamical picture for 
$g\in\mathcal{U}_+$:   

\begin{figure}[hbt!]
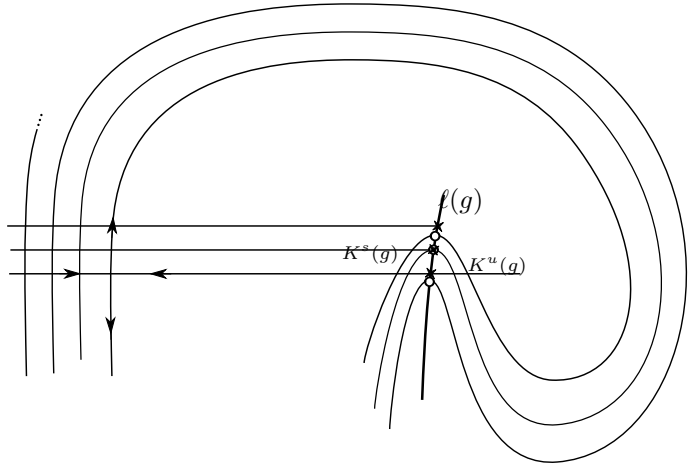
\caption{The line of tangencies $\ell(g)$ and the Cantor sets $K^s(g)$ and $K^u(g)$ for $g\in\mathcal{U}_+$: the crosses are points in $K^s(g)$ and the dots are points in 
$K^u(g)$.}\label{f.6}
\end{figure}

In particular, the problem of \emph{persistent} tangencies for \emph{all} $g\in\mathcal{U}_+$, i.e., the issue that the stable and unstable laminations of $K_g$ meet tangentially at some point in $V$ for all $g\in\mathcal{U}_+$, is reduced to the question of studying sufficient conditions for a non-trivial intersection $K^s(g)\cap K^u(g)\neq\emptyset$ between the Cantor sets $K^s(g)$ and $K^u(g)$ of the line $\ell(g)$. 

\subsubsection{Intersections of thick Cantor sets and Newhouse gap lemma}\label{ss.gap-lemma}

By thinking of $\ell(g)$ as a subset of the real line $\mathbb{R}$, our current objective is clearly equivalent to producing sufficient conditions to ensure that two Cantor sets in $\mathbb{R}$ have non-trivial intersection. Keeping this goal in mind, S. Newhouse introduced a notion of \emph{thickness} $\tau(C)$ of a Cantor set 
$K\subset\mathbb{R}$: 

\begin{definition}Let $K\subset\mathbb{R}$ be a Cantor set. A \emph{gap} $U$ of $K$ is a connected component of $\mathbb{R}-K$, and a \emph{bounded gap} $U$ of $K$ is a bounded connected component of $\mathbb{R}-K$. Given $U$ a bounded gap of $K$ and $u\in\partial U$, the \emph{bridge} $C$ of $K$ at $u$ is the maximal interval $C\subset\mathbb{R}$ such that $u\in\partial C$ and $C$ contains \emph{no} gap $U'$ of $K$ with length $|U'|$ greater than or equal to the lenght $|U|$ of $U$. In this language, the thickness of $K$ at $u$ is $\tau(K,u):=|C|/|U|$ and the \emph{thickness} $\tau(K)$ of $K$ is 
$$\tau(K)=\inf\limits_{u}\tau(K,u).$$
\end{definition}
The thickness is a nice concept for our purposes because of the following fundamental result nowadays called \emph{gap lemma}:

\begin{lemma}[Newhouse gap lemma] Let $K$ and $K'$ be two thick Cantor sets of $\mathbb{R}$ in the sense that $\tau(K)\cdot\tau(K')>1$. Then, one of the following possibilities occur:
\begin{itemize}
\item $K'$ is contained in a gap of  $K$ (i.e., a connected component of $\mathbb{R}-K$),
\item $K$ is contained in a gap of $K'$,
\item $K\cap K'\neq\emptyset$.
\end{itemize}
\end{lemma}

\begin{remark}\label{r.linked-Cantors} A practical way of using Newhouse gap lemma by looking at the relative position of two Cantor sets is the following. We say that two Cantor sets $K,K'\subset\mathbb{R}$ are \emph{linked} whenever their convex hulls $I, I'$ are linked in the sense that $I\cap I'\neq\emptyset$ but neither $I\not\subset I'$ nor $I'\not\subset I$. Then, by Newhouse's gap lemma, two linked Cantor sets $K,K'\subset\mathbb{R}$ with $\tau(K)\cdot\tau(K')>1$ must intersect non-trivially because, as the reader can easily check, the assumption that $K$ and $K'$ are linked prohibits a gap of $K'$ to contain $K$ and vice-versa. 
\end{remark}

The proof of this lemma is not difficult and one can find it on page 63 of \cite{PT} for instance. Of course, the gap lemma put us in position to come back to the discussion of persistence of tangencies for $g\in\mathcal{U}_+$. Indeed, the statement of the gap lemma hints that one has persistence of tangencies for all $g\in\mathcal{U}_+$ as soon as $\tau(K^s)\cdot \tau(K^u)>1$ for the initial dynamics $f\in\mathcal{U}_0$. Actually, this fact was shown to be true by S. Newhouse, but this is not an immediate consequence of his gap lemma because we need to know that the Cantor sets $K^s(g)$ and $K^u(g)$ are thick for all $g\in\mathcal{U}_+$ and we have only that $K^s$ and $K^u$ are thick.

At this point, the idea (already mentioned above) is to play with continuity of dynamical objects: intuitively, $K^s(g)$ and $K^u(g)$ must be thick because they are close to the thick Cantor sets $K^s$ and $K^u$. However, the formal implementation of this idea is rather technical and we'll content ourselves with a mere description of the crucial points of the argument. 

\subsubsection{Continuity of thickness and Newhouse phenomena}

Firstly, one has to explain what does it mean for $\mathcal{F}^s(g)$ to be ``close'' to $\mathcal{F}^s(f)$: for our purposes, we'll say that $\mathcal{F}^s(g)$ is close to $\mathcal{F}^s(f)$ whenever the map $U\ni x\mapsto T_x\mathcal{F}^s(f)(x)$ is $C^1$-close to the map $U\ni x\mapsto T_x\mathcal{F}^s(g)(x)$. Here, $U$ is our preferred (small) neighborhood of the horseshoe $K$ of $f$ and $\mathcal{F}^s(f)(x)$, resp. $\mathcal{F}^s(g)(x)$, is the leaf of $\mathcal{F}^s(f)$, resp. $\mathcal{F}^s(g)$ passing through $x$. 

As it is shown in Theorem 8 of Appendix 1 of \cite{PT}, $\mathcal{F}^s(g)$ is close to $\mathcal{F}^s(f)$ when $g$ is $C^k$-close to $f$ for $k\geq 2$. Using this result, one can show that the line of tangency $\ell(g)$ is $C^1$-close to $\ell(f)$, and the projections $\pi^s_g$ and $\pi^u_g$ are $C^1$-close to $\pi^s_f$ and $\pi^u_f$. An immediate consequence of this is: one can nicely identify both $\ell(g)$ and $\ell$ with the interval $I=[0,1]\subset\mathbb{R}$ in such a way that $K^s(g)$ is close to $K^s (=K^s(f))$ in the \emph{Hausdorff topology}\footnote{We say that a compact subset $A\subset\mathbb{R}$ is $\delta$-close to a compact subset $B\subset\mathbb{R}$ in the Hausdorff topology when for each $y\in B$ there exists $x\in A$ with $|y-x|<\delta$, and for each $w\in A$ there exists $z\in B$ with $|w-z|<\delta$.}. However, this is not very useful because it is not true in general that $\tau(K^s(g))$ is close to $\tau(K^s)$ when $K^s(g)$ and $K^s$ are close in Hausdorff topology. Here, the idea to overcome this difficulty relies on the observation that $K^s(g)$ and $K^s$ belong to a special class of Cantor sets known as \emph{regular (dynamically defined) Cantor sets}. 

More precisely, we say that a Cantor set $K\subset\mathbb{R}$ is $C^r$-\emph{regular}, $r\geq 1$, if there are disjoint compact intervals $I_1,\dots,I_l\subset\mathbb{R}$ and a uniformly expanding $C^r$ function $\psi:I_1\cup\dots\cup I_l\to I$ (i.e., $|\psi'(x)|>1$ for any $x$) from the disjoint union $I_1\cup\dots\cup I_l$ to its convex hull $I$ such that: 
\begin{itemize}
\item $K=\bigcap\limits_{n\in\mathbb{N}}\psi^{-n}(I)$, that is, $K$ is defined by the dynamics $\psi$, and 
\item the collection $\{I_1,\dots,I_l\}$ is a \emph{Markov partition}: for any $1\leq j\leq l$, the interval $\psi(I_j)$ is the convex hull of the union of some of the intervals $I_i$ and 
$\psi^{n(j)}(I_j)\supset I_1\cup\dots\cup I_l$ for some large $n(j)\in\mathbb{N}$.
\end{itemize}

Regular (dynamically defined) Cantor sets are very common in nature. For example, the classical \emph{ternary Cantor set} $K_0$ is a regular Cantor set. Indeed, it is not hard to see that $K_0=\bigcap\limits_{n\in\mathbb{N}}\psi^{-n}([0,1])$ where $\psi:[0,1/3]\cup [2/3,1]\to [0,1]$ is the (piecewise affine) expanding function defined by
$$\psi(x)=\left\{\begin{array}{ll}3x & \textrm{ if } x\in [0,1/3] \\ 3x-2 & \textrm{ if } x\in [2/3,1] \end{array}\right.$$
Other examples are the Cantor sets $W^s(p)\cap K$ and $W^s(p_g)\cap K_g$: as it is shown in Chapter 4 of \cite{PT}, they are $C^k$-regular Cantor sets whenever $f$ and $g$ are $C^k$-diffeomorphisms.

The class of $C^r$-regular Cantor sets admit a natural $C^r$-\emph{topology}: we say that two $C^r$-regular Cantor sets $K$ and $\widetilde{K}$ are $C^r$-\emph{close} whenever the extremal points of the associated intervals $I_1,\dots, I_l$ and $\widetilde{I_1},\dots,\widetilde{I_l}$ are close and the expanding functions $\psi$ and $\widetilde{\psi}$ are $C^r$-close. For example, the $C^k$-regular Cantor sets $W^s(p)\cap K$ and $W^s(p_g)\cap K_g$ are $C^k$-close when $f$ and $g$ are $C^k$-close.

These definitions are well-adapted to the study of homoclinic tangencies because of the following fundamental fact:

\begin{proposition}\label{p.thickness-continuity} The thickness of $C^k$-regular Cantor sets vary \emph{continuously} in the $C^k$-topology for $k\geq 2$.
\end{proposition} 
See Chapter 4 of \cite{PT} for more discussion on this proposition.

Therefore, if $f\in\mathcal{U}_0$ is a $C^k$-diffeomorphism, $k\geq 2$, and $\tau(W^s_{loc}(p)\cap K)\cdot\tau(W^u_{loc}(p)\cap K)>1$, then $\tau(W^s_{loc}(p_g)\cap K_g)\cdot\tau(W^u_{loc}(p_g)\cap K_g)>1$ for all $g\in\mathcal{U}_+$. Thus, since $\ell(g)$ is $C^1$-close to $\ell(f)$, and $\pi^s(g)$ and $\pi^u(g)$ are $C^1$ close to $\pi^s(f)$ and 
$\pi^u(f)$, it is possible to compare $W^s(p_g)\cap K_g$ and $K^s(g)$ via a $C^1$-diffeomorphism whose derivative is $C^0$-close to the identity. In particular, since $C^1$-diffeomorphisms with a derivative $C^0$-close to the identity don't change in a drastic way the thickness, one has that $\tau(K^s(g))$ is close $\tau(W^s(p_g)\cap K_g)$, so that 
$\tau(K^s(g))\cdot \tau(K^u(g))$ is close to $\tau(W^s_{loc}(p_g)\cap K_g)\cdot\tau(W^u_{loc}(p_g)\cap K_g)>1$. Hence, we conclude that $\tau(K^s(g))\cdot \tau(K^u(g))>1$ for all $g\in\mathcal{U_+}$. Moreover, for each $g\in\mathcal{U}_+$, the Cantor sets $K^s(g)$ and $K^u(g)$ of the line $\ell(g)$ are linked (as one can see from Figure \ref{f.6}). Thus, by Newhouse's gap lemma (or more precisely Remark \ref{r.linked-Cantors}), we obtain that $K^s(g)\cap K^u(g)\neq\emptyset$ for all $g\in\mathcal{U}_+$. 

In other words, we have just outlined the proof of the following result about persistance of tangencies:

\begin{theorem}\label{t.persistent-tangency} Let $f\in\mathcal{U}_0$ be a $C^k$-diffeomorphism, $k\geq 2$, and suppose that $\tau(W^s_{loc}(p)\cap K)\cdot\tau(W^u_{loc}(p)\cap K)>1$. Then, for all $g\in\mathcal{U}_+$, the stable and unstable laminations of the horseshoe $K_g$ intersect tangentially at some point in $V$.
\end{theorem} 

Once we dispose of this theorem on persistence of tangencies in our toolbox, we're ready to discuss the Newhouse phenomena. Again, we start with a $C^k$-diffeomorphism $f\in\mathcal{U}_0$ with $k\geq 2$ and we now \emph{assume} that:
\begin{itemize}
\item the periodic point $p$ is \emph{dissipative}, i.e., $|\det df^n(p)|\neq 1$ where $n$ is the period of $p$, and
\item $\tau(W^s_{loc}(p)\cap K)\cdot\tau(W^u_{loc}(p)\cap K)>1$
\end{itemize} 
For sake of concreteness, let's suppose that $|\det df^n(p)|< 1$. Note that this implies $|\det dg^n(p_g)|< 1$ for all $g\in\mathcal{U}$. Denote by $\lambda_g<1<\sigma_g$ the eigenvalues of $dg^n(p)$, so that $|\lambda_g\cdot\sigma_g|=|\det dg^n(p_g)|<1$. Now, given $g_0\in\mathcal{U}_+$, we know by Theorem \ref{t.persistent-tangency} that the stable and unstable laminations of $K_{g_0}$ meet tangentially at some point in $V$. Since the stable and unstable manifolds of the periodic point $p_{g_0}$ are dense in the stable and unstable laminations of $K_{g_0}$ (cf. Chapter 2 of \cite{PT}), one can apply arbitrarily small perturbations to $g_0\in\mathcal{U}_+$ so that there is no loss of generality in assuming that $W^s(p_{g_0})$ and $W^u(p_{g_0})$ meet tangentially at some point in the region $V$. Starting from this quadratic tangency, one gets the following picture for diffeomorphisms $g_{\mu}\in\mathcal{U}_+$ close to $g_0$: 

\begin{figure}[hbt!]
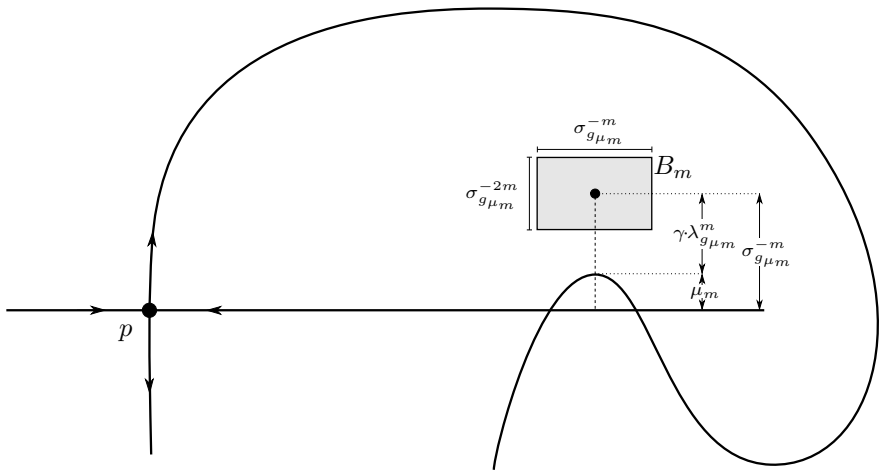
\caption{Selection of adequate boxes $B_m$ for the renormalization of $g_{\mu_m}$.}\label{f.7}
\end{figure}

As it is shown in Chapter 3 of \cite{PT}, one can carefully choose parameters\footnote{In principle, the parameters $\mu$ must vary in some infinite-dimensional manifold in order to $g_{\mu}$ parametrize a neighborhood of $g_0$, but for sake of simplicity of the exposition, we will think of this parameter as a \emph{real number} $\mu\in\mathbb{R}$ measuring the distance between the line $W^s(p_{g_\mu})\cap V$ and the tip of the parabola $W^u(p_{g_{\mu}})\cap V$ as indicated in Figure \ref{f.7}.} $\mu_m$ ($m\in\mathbb{N}$) such that 
\begin{itemize}
\item $\mu_m\to 0$ as $m\to+\infty$ and 
\item the map $g_{\mu_m}^m|_{B_m}$ can be \emph{renormalized}\footnote{That is, one can perform an adequate $\mu_m$-\emph{dependent} change of coordinates 
$\phi_{\mu_m}$ on $g_{\mu}^m|_{B_m}$ to get a new dynamics $G_m=\phi_{\mu_m}^{-1}\circ g_{\mu}^m|_{B_m}\circ\phi_{\mu_m}$.} in such a way that the renormalizations $G_m$ of $g_{\mu_m}^m|_{B_m}$ $C^2$-converge to the endomorphism $(\widetilde{x},\widetilde{y})\mapsto (\widetilde{y},\widetilde{y}^2)$. 
\end{itemize}
Note that the diffeomorphisms $G_m$ are converging to an endomorphism and this may seem strange at first sight. However, this is natural in view of the assumption that the periodic point $p$ is dissipative: in fact, the area-contraction condition $|\det df^n(p)|<1$ says that $g_{\mu_m}^m$ become strongly area-contracting as $m\to\infty$ and consequently $g_{\mu_m}^m|_{B_m}$ converges to a curve and $G_m$ converges to an endomorphism of this curve as $m\to\infty$. 

Next, we observe that the endomorphism $(\widetilde{x},\widetilde{y})\mapsto (\widetilde{y},\widetilde{y}^2)$ has an attracting fixed point at $(\widetilde{x},\widetilde{y})=(0,0)$. Therefore, by $C^2$ convergence of $G_m$ towards this endomorphism, we conclude that $g_{\mu_m}^m$ has an attracting fixed point in $B_m$ for all $m$ sufficiently large. In other words, $g_{\mu_m}$ has a \emph{sink} (attracting periodic point) in the region $V$ for all $\mu_m$ sufficiently small. 

This last statement can be reformulated as follows. For each $m\in\mathbb{N}$, denote by $R_m=\{g\in\mathcal{U}_+: g \textrm{ has } m \textrm{ sinks}\}$. Note that $R_m$ is open for all $m\in\mathbb{N}$ (because any sink is persistence under small perturbations of the dynamics). Moreover, since $g_0\in\mathcal{U}_+$ was arbitrary in the previous argument, we also have that $R_1$ is dense in $\mathcal{U}_+$.

At this stage, the idea of S. Newhouse is to iterate this argument to show that the set 
$$R_{\infty}:=\bigcap\limits_{m\in\mathbb{N}} R_m$$
of diffeomorphisms of $\mathcal{U}_+$ with \emph{infinitely many sinks} is \emph{residual}\footnote{A set is residual in the Baire category sense whenever it contains a countable intersection of open and dense subsets.} in Baire sense (and, in particular, $R_{\infty}$ is dense in $\mathcal{U}_+$). Since $R_m$ is open in $\mathcal{U}_+$ for all $m\in\mathbb{N}$ and $R_1$ is dense in $\mathcal{U}_+$, it suffices to prove that $R_{m+1}$ is dense in $R_m$ for all $m\in\mathbb{N}$ to conclude that $R_{\infty}$ is residual. 

In this direction, one starts with $g_0\in R_m$ with $m$ periodic sinks $\mathcal{O}_1(g_0),\dots,\mathcal{O}_m(g_0)$. By Theorem \ref{t.persistent-tangency}, we know that the stable and unstable laminations of $K_{g_0}$ meet tangentially somewhere in $V$. Since $W^s(p_{g_0})$, resp. $W^u(p_{g_0})$, is dense in the stable, resp. unstable, lamination of $K_{g_0}$, we can assume (up to performing an arbitrarily small perturbation on $g_0$) that $W^s(p_g)$ and $W^u(p_g)$ meet tangentially at some point $q_{g_0}\in V$ and $g_0$ has $m$ periodic sinks. Next, we select $T$ a small neighborhood of $q_{g_0}$ such that none of the periodic sinks passes through $W$, i.e., $W\cap \mathcal{O}_i(g_0)=\emptyset$ for each $i=1,\dots,m$. By repeating the ``renormalization'' arguments above (with $V$ replaced by $T$), one can produce a sequence of diffeomorphisms $(g_{\mu_j})_{j\in\mathbb{N}}$ converging to $g_0$ as $j\to\infty$ such that $g_{\mu_j}$ has a sink $\mathcal{O}(g_{\mu_j})$ passing through $T$. Because the sinks $\mathcal{O}_i(g_{\mu_j})$ don't pass through $T$ for all $j$ sufficiently large, this means that $\mathcal{O}(g_{\mu_j})$ is a new sink of $g_{\mu_j}$, that is, we obtain that $g_{\mu_j}\in R_{m+1}$ for all $j$ sufficiently large. Since $g_{\mu_j}\to g_0$ as $j\to\infty$, we conclude that $R_{m+1}$ is dense in $R_m$.

Thus, we have given a sketch of the proof of the following result:

\begin{theorem}[S. Newhouse]\label{t.Newhouse} Let $k\in\mathbb{N}$, $k\geq 2$, and let $f\in\mathcal{U}_0$ be a $C^k$-diffeomorphism such that 
\begin{itemize}
\item the periodic point $p$ is \emph{dissipative}, say $|\det df^n(p)|< 1$ where $n$ is the period of $p$, and
\item $\tau(W^s_{loc}(p)\cap K)\cdot\tau(W^u_{loc}(p)\cap K)>1$
\end{itemize}
Then, the subset $R_{\infty}\subset\mathcal{U}_+$ of diffeomorphisms with infinitely many sinks is residual. 
\end{theorem}

This result of coexistence of infinitely many sinks for a residual (and, hence, dense) subset of diffeomorphisms of $\mathcal{U}_+$ is the so-called \emph{Newhouse phenomena}. This theorem is very important because it says that for a topologically big (residual) set $R_{\infty}$ of diffeomorphisms the dynamics is so complicated that there are infinitely many attractors. Thus, if we pick at random a point $x$ of $U\cup V$, it is very hard to decide (from the computational point of view for instance) the future of the orbit of $x$ because it can be attracted by any one of the infinitely many sinks. 

In other words, the Newhouse phenomena says that it is not reasonable to try to understand the local dynamics of \emph{all} $g\in\mathcal{U}_+$. At this point, since we know that it is too naive to try to dynamically describe all $g\in\mathcal{U}_+$, we can ask: 
\begin{center}
\textit{What about the local dynamical behavior of \emph{most} $g\in\mathcal{U}_+$?}
\end{center} 

The discussion of this question will occupy the remainder of this text. For now, we close this subsection with two comments:

\begin{remark}\label{r.Newhouse-w/o-thickness} Actually, S. Newhouse proved in \cite{N3} (see also Chapter 6 of \cite{PT}) that one can remove the second assumption (on thicknesses) in the statement of Theorem \ref{t.Newhouse}: indeed, starting with any dissipative area-contracting hyperbolic periodic point $p$ (of saddle-type) of a $C^2$-diffeomorphism $f$ having some point $q$ of tangency between $W^s(p)$ and $W^u(p)$, S. Newhouse can construct open sets $\mathcal{U}$ arbitrarily close to $f$ such that the subset of diffeomorphisms of 
$\mathcal{U}$ with infinitely many sinks is residual in $\mathcal{U}$.
\end{remark} 

\begin{remark}The attentive reader certainly noticed that we insisted that Newhouse phenomena (Theorem \ref{t.Newhouse}) concerns $C^2$-diffeomorphisms. In fact, this regularity assumption is crucial to get the continuity of the thickness of regular Cantor sets in Proposition \ref{p.thickness-continuity}. Indeed, the proof of this proposition in Chapter 4 of \cite{PT} strongly relies on the so-called \emph{bounded distortion} property saying that the shape of gaps and bridges of a $C^2$-regular Cantor set is ``essentially constant in all scales''. Of course, the continuity of thickness is one of the central mechanisms for Newhouse phenomena (as it ensures that the Cantor sets $K^s(g)$ and $K^u(g)$ intersect for all $g\in\mathcal{U}_+$) and the reader may be curious whether the Newhouse phenomena occurs for $C^1$-diffeomorphisms. As a matter of fact, it is known that the thickness of $C^1$-regular Cantor sets is not continuous, so that the Newhouse gap lemma can not be applied in the $C^1$-context. Of course, it could be that $C^1$-regular Cantor sets intersect often in a stable manner, thus giving some hope for an analog of Newhouse's thickness mechanism to survive in the $C^1$-setting. However, this possibility was recently dismissed by C. (Gugu) Moreira \cite{Mo} and this ``absence of Newhouse mechanism'' was used by C. Moreira, E. Pujals and the author \cite{MMP} to check that among certain families of dynamical systems it is possible to get a sort of Newhouse phenomena in the $C^2$-setting but still the $C^1$-generic element of the family has finitely many sinks.
\end{remark}

\subsection{Homoclinic bifurcations associated to thin horseshoes}\label{ss.Newhouse-Palis-Takens} Let's come back to the setting of the beginning of Subsection \ref{ss.bifurcations-Newhouse}, that is, $f:M\to M$ is a $C^2$-diffeomorphism of a surface $M$ with a horseshoe $K$ and a periodic point $p\in K$ whose stable and unstable manifolds have a quadratic tangency at some point $q\in M-K$. Consider again $U$ a sufficiently small neighborhood of $K$ and $V$ a sufficiently small neighborhood of the orbit of $q$, and let's fix $\mathcal{U}$ a sufficiently small $C^2$-neighborhood of $f$ organized into the open sets $\mathcal{U}_-$ and $\mathcal{U}_+$ and the codimension $1$ hypersurface $\mathcal{U}_0$ depending on the relative positions of $W^s(p_g)$ and $W^u(p_g)$ near $V$. 

From now on, we will be interested in the local dynamics of $\Lambda_g$ for most $g\in\mathcal{U}_+$. Of course, there are plenty of reasonable ways of formalizing the notion of ``most'' here. For the sake of these notes, we will adopt the following definition:

\begin{definition}We say that a subset (i.e., a property) $\mathcal{P}\subset\mathcal{U}_+$ contains (i.e., holds for) \emph{most} $g\in\mathcal{U}_+$ whenever for every smooth $1$-parameter family $(g_t)_{|t|<t_0}$ with 
\begin{itemize}
\item $g_t\in\mathcal{U}_-$ for $-t_0<t<0$, $g_0\in\mathcal{U}_0$ and $g_t\in\mathcal{U}_+$ for $0<t<t_0$, and
\item $(g_t)_{|t|<t_0}$ is transverse to the codimension $1$ hypersurface $\mathcal{U}_0$
\end{itemize}
one has that\footnote{Here, $\textrm{Leb}$ is the $1$-dimensional Lebesgue measure.} 
$$\lim\limits_{t\to 0^+}\frac{\textrm{Leb}(\{s\in(0,t): g_s\in\mathcal{P}\})}{t}=1$$
\end{definition}

In plain terms, $\mathcal{P}\subset \mathcal{U}_+$ contains most $g\in\mathcal{U}_+$ if $\mathcal{P}$ has density $1$ at $\mathcal{U}_0$, where the density is measured along smooth generic (i.e., transverse) $1$-parameter families crossing $\mathcal{U}_0$. 

Using this reasonable notion of ``most $g\in\mathcal{U}_+$'', the following question makes sense:
\begin{center}
\textit{Is $\Lambda_g$ a (uniformly hyperbolic) horseshoe for most $g\in\mathcal{U}_+$?}
\end{center}

From our previous experience with the Newhouse phenomena, we know that this question is delicate when the horseshoe $K$ is thick: indeed, we saw that if the Cantor sets $W^s(p)\cap K$ and $W^u(p)\cap K$ are thick, then one has persistence of tangencies in $\mathcal{U}_+$ and this is a dangerous scenario conspiring against the hyperbolicity of $\Lambda_g$. On the other hand, it is intuitive that if the horseshoe $K$ is ``thin'' in some adequate sense, one can get rid of tangencies and this gives some hope that in this situation $\Lambda_g$ is a horseshoe for most $g\in\mathcal{U}_+$. 

In fact, the intuition of the previous paragraph can be formalized with the aid of the notion of \emph{Hausdorff dimension}.

\begin{definition} Let $A\subset\mathbb{R}^n$. Given $\mathcal{O}=\{O_i\}_{i\in I}$ a countable open cover\footnote{I.e., the subsets $O_i$ are open and $A\subset\bigcup\limits_{i\in I}O_i$.} of $A$, we define its diameter $\textrm{diam}(\mathcal{O})=\sup\limits_{i\in I}\textrm{diam}(O_i)$. Given a real number $\alpha\geq 0$, the \emph{Hausdorff $\alpha$-measure} of $A$ is 
$$m_{\alpha}(A):=\lim\limits_{\varepsilon\to 0} \,\, \inf\limits_{\substack{\mathcal{O} \textrm{ countable open cover of } A \\ \textrm{with } \textrm{diam}(\mathcal{O})<\varepsilon}} \, \, \sum\limits_{O_i\in\mathcal{O}} \textrm{diam}(O_i)^{\alpha}$$
and the \emph{Hausdorff dimension} of $A$ is $\textrm{HD}(A):=\inf\{\alpha\geq 0: m_{\alpha}(A)=0\}$.
\end{definition}

As an exercise the reader can try to show from the definitions that the Hausdorff dimension has the following general properties:

\begin{proposition}\label{p.Hausdorff-dimension} The Hausdorff dimension has the following properties:
\begin{itemize}
\item[(a)] it is monotone: $\textrm{HD}(B)\leq \textrm{HD}(A)$ whenever $B\subset A$; 
\item[(b)] it is countably stable: $\textrm{HD}(\bigcup\limits_{i\in\mathbb{N}} A_i)=\sup\limits_{i\in\mathbb{N}}\textrm{HD}(A_i)$;
\item[(c)] $\textrm{HD}(A)=0$ whenever $A\subset\mathbb{R}^n$ is finite or countable;
\item[(d)] a compact $m$-dimensional submanifold $M^m\subset\mathbb{R}^n$ has Hausdorff dimension $\textrm{HD}(M^m)=m$; 
\item[(e)] it doesn't increase under Lipschitz maps, i.e., $\textrm{HD}(f(A))\leq \textrm{HD}(A)$ if $f:A\to\mathbb{R}^k$ is a Lipschitz\footnote{I.e., there exists a constant $L>0$ such that $|f(x)-f(y)|\leq L|x-y|$ for all $x,y$.} map; 
\item[(f)] the Hausdorff dimension of a product set $A\times B$ satisfies $\textrm{HD}(A\times B)\geq \textrm{HD}(A)+\textrm{HD}(B)$;
\item[(g)] any measurable $A\subset\mathbb{R}^n$ with $\textrm{HD}(A)<n$ has zero Lebesgue measure.
\end{itemize}
\end{proposition}

Coming back to the study of the local dynamics of $C^2$-diffeomorphisms $g\in\mathcal{U}$, let's consider the horseshoe $K_g$. We define the \emph{stable} (resp. \emph{unstable}) \emph{dimension}\footnote{We measure the stable dimension of $K_g$ using the unstable manifold of $p_g$ because we're interested in the transverse structure of the stable set of $K_g$. Also, we call $d_s(K_g)$ the stable dimension of $K_g$ instead of stable dimension of $K_g$ at $p_g$ because it is possible to prove that $d_s(K_g)=W^u(x)\cap K_g$ for all $x\in K_g$.} $d_s(K_g)$ (resp. $d_u(K_g)$) of $K_g$ as the Hausdorff dimension of $W^u(p_g)\cap K_g$ (resp. $W^s(p_g)\cap K_g$). Also, for later use, we denote by $d_s^0$ and $d_u^0$ the stable and unstable Hausdorff dimensions of $K=K_f$, $f\in\mathcal{U}_0$. 

The stable and unstable dimensions $d_s(K_g)$ and $d_u(K_g)$ are nicely related to the geometry of the horseshoe $K_g$ because of the formula\footnote{The idea behind this formula is that at small scales the horseshoe $K_g$ looks like the product of the regular Cantor sets $W^s_{loc}(x)\cap K_g$ and $W^u_{loc}(x)\cap K_g$ and this allows to get an improved version of item (f) of Proposition \ref{p.Hausdorff-dimension} in this case.}:
\begin{equation}\label{e.Hausdorff-horseshoe}
\textrm{HD}(K_g)=d_s(K_g)+d_u(K_g)
\end{equation}

Using the notion of stable and unstable dimensions, S. Newhouse, J. Palis and F. Takens \cite{NP}, \cite{PT2} in 1987 proved the following theorem:

\begin{theorem}[S. Newhouse, J. Palis and F. Takens]\label{t.Newhouse-Palis-Takens} Suppose that $d_s^0+d_u^0<1$ for $f\in\mathcal{U}_0$. Then, $\Lambda_g$ is a uniformly hyperbolic horseshoe for most $g\in\mathcal{U}_+$.
\end{theorem}

Informally speaking, this theorem says that if the horseshoe $K=K_f$ is \emph{thin} enough in the sense that $\textrm{HD}(K)=d_s^0+d_u^0<1$, then $\Lambda_g$ is a horseshoe for most ways of unfolding the quadratic tangency of $f$ at $q$ (i.e., for most $g\in\mathcal{U}_+$). 

Let us now explain why this theorem is intuitively plausible. Let us consider $(g_t)_{|t|<t_0}$ a smooth $1$-parameter family transverse to $\mathcal{U}_0$ at $f=g_0$. Our first obstacle towards hyperbolicity is the issue of tangencies. So, using the notations from Sub-subsection \ref{sss.tangencies-Cantors}, let us again consider the regular Cantor sets $K^s(g_t)$ and $K^u(g_t)$ on the line of tangencies $\ell(g_t)$ whose intersections $K^s(g_t)\cap K^u(g_t)$ account for all tangencies between the stable and unstable laminations of $K_{g_t}$. Because the tangencies for $f\in\mathcal{U}_0$ are quadratic, by adequate reparametrization, we may think that the regular Cantor sets $K^s(g_t)$ and $K^u(g_t)$ live in the real line $\mathbb{R}$ and they \emph{move} with \emph{unit speed} relatively to each other, i.e., 
$$K^s(g_t)=K^s(g_0) \quad \textrm{ and } \quad K^u(g_t)=K^u(g_0)+t$$
for all $|t|<t_0$. In this context, note that $K^s(g_t)\cap K^u(g_t)\neq\emptyset$ if and only if $t\in K^s(g_0)\ominus K^u(g_0)$ where $A\ominus B:=\{x-y: x\in A,\, y\in B\}$ denotes the arithmetic difference between $A\subset\mathbb{R}$ and $B\subset\mathbb{R}$. In other words, the arithmetic difference $K^s(g_0)\ominus K^u(g_0)$ of the regular Cantor sets $K^s(g_0)$ and $K^u(g_0)$ accounts for all parameters $t\in(-t_0,t_0)$ such that the stable and unstable laminations of $K_{g_t}$ exhibits some tangency. Therefore, it is desirable to know the size of this arithmetic difference. 

In this direction, one observes that $K^s(g_0)\ominus K^u(g_0)=\pi(K^s(g_0)\times K^u(g_0))$ where $\pi:\mathbb{R}^2\to\mathbb{R}$ is the projection $\pi(x,y)=x-y$. Since $K^s(g_0)$ and $K^u(g_0)$ are regular Cantor sets of Hausdorff dimensions\footnote{This holds because the regular Cantor sets $W^s_{loc}(p)\cap K$ and $W^u_{loc}(p)\cap K$ have Hausdorff dimension $d_s^0$ and $d_u^0$ by definition, and $K^s(g_0)$ and $K^u(g_0)$ are diffeomorphic to  $W^s_{loc}(p)\cap K$ and $W^u_{loc}(p)\cap K$ (so that item (e) of Proposition \ref{p.Hausdorff-dimension} applies).} $d_s^0$ and $d_u^0$, one has that the product set $K^s(g_0)\times K^u(g_0)$ has Hausdorff dimension $d_s^0+d_u^0$. By item (e) of Proposition \ref{p.Hausdorff-dimension}, we obtain that the arithmetic difference $K^s(g_0)\ominus K^u(g_0)=\pi(K^s(g_0)\times K^u(g_0))$ has Hausdorff dimension 
$$\textrm{HD}(K^s(g_0)\ominus K^u(g_0))\leq d_s^0+d_u^0$$
because the projection $\pi$ is Lipschitz. By item (g) of Proposition \ref{p.Hausdorff-dimension}, we conclude that the arithmetic difference $K^s(g_0)\ominus K^u(g_0)\subset\mathbb{R}$ has zero Lebesgue measure. In other words, the assumption $d_s^0+d_u^0<1$ imposes a severe restriction on the set of parameters $|t|<t_0$ such that the invariant laminations of $K_{g_t}$ exhibits a tangency, namely, these parameters have zero Lebesgue measure. 

At this point, we got rid of the issue of tangencies (from the measure-theoretical point of view), but unfortunately this is not sufficient to ensure the hyperbolicity of $\Lambda_g$: indeed, while it is quite clear that the pieces of orbits passing near the horseshoe $K_g$ have natural canditates for the stable and unstable directions $E^s(x)$ and $E^u(x)$ in Definition \ref{d.hyp-set} (of hyperbolicity), this is not so clear in the region $V$ (near the quadratic tangency for $f\in\mathcal{U}_0$) as the candidate directions $E^s(x)$ and $E^u(x)$ may reverse their role (and thus the hyperbolicity is lost) due to almost tangencies between the invariant laminations of $K_g$. In other words, we need not only to ensure that $K^s(g_t)\cap K^u(g_t)=\emptyset$, but we also have to ensure that $K^s(g_t)$ and $K^u(g_t)$ are sufficiently far apart from each other in order to obtain the hyperbolicity of $\Lambda_{g_t}$. 

To formalize the idea of the previous paragraph, one needs to know the localization of points of $\Lambda_{g_t}-K_{g_t}$, that is, the points of $\Lambda_{g_t}$ whose $g_t$-orbit passes by the region $V$. Here, one can show (see the proposition by the end of page 213 of \cite{PT}) that, given $c>0$, all points of $\Lambda_{g_t}-K_{g_t}$ have some $g_t$-iterate in $V$ at a distance $\leq c\cdot t$ of the invariant laminations of $K_{g_t}$ for all $t$ sufficiently small (depending on $c>0$). This fact is very interesting because it says that one can understand the orbits in $\Lambda_{g_t}-K_{g_t}$ by looking at the intersection of the $c\cdot t$-neighborhoods $\mathcal{F}^s(t)$ and $\mathcal{F}^u(t)$ of the stable and unstable laminations of $K_{g_t}$. We illustrate this intersection in the figure below:

\begin{figure}[hbt!]
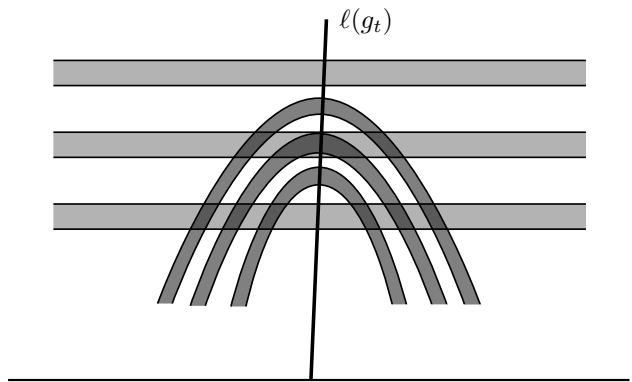
\caption{Neighborhoods $\mathcal{F}^s(t)$/$\mathcal{F}^u(t)$ indicated by horizontal/parabolic-like strips. The intersection $\mathcal{F}^s(t)\cap\mathcal{F}^u(t)$ is indicated by strong grey.}\label{f.8}
\end{figure}

In this picture, we are again using the fact that the tangencies for $f\in\mathcal{U}_0$ are quadratic. From this figure, we see that the geometry of $\mathcal{F}^s(t)\cap\mathcal{F}^u(t)$ is controlled by the relative position of the $ct$-neighborhoods $A_{ct}$ and $B_{ct}$ of the Cantor sets $K^s(g_t)$ and $K^u(g_t)$ on the line of tangency $\ell(g_t)$: for instance, \emph{if} the distance between $A_{ct}$ and $B_{ct}$ is $\geq 2ct$, then the angle between the leaves of the stable and unstable foliations at any point $x\in \mathcal{F}^s(t)\cap\mathcal{F}^u(t)$ is $\geq \sqrt{2ct}$. See the figure below.

\begin{figure}[hbt!]
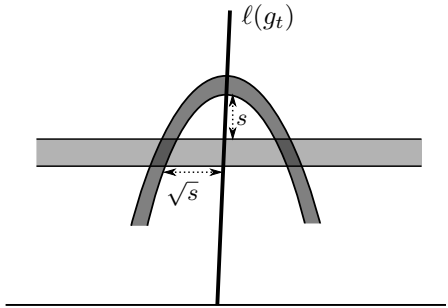
\caption{Angle estimate in the region $\mathcal{F}^s(t)\cap\mathcal{F}^u(t)$ (here $s:=2ct$).}\label{f.9}
\end{figure}

In other terms, if the distance between the sets $A_{ct}$ and $B_{ct}$ is $\geq 2ct$, then we don't see almost tangencies, i.e., the angle between leaves of the stable and unstable foliations is $\geq \sqrt{2ct}$. Since the tangents to the leaves of the stable and unstable foliations at $\Lambda_g-K_g$ are the natural candidates for stable and unstable directions over $\Lambda_{g_t}$ in the sense of Definition \ref{d.hyp-set}, it is not surprising that one can actually prove that $\Lambda_{g_t}$ is a hyperbolic set (and hence a horseshoe) whenever this angle estimate holds. 

Therefore, we have that if the distance $d(A_{ct},B_{ct})$ between $A_{ct}$ and $B_{ct}$ is $\geq 2ct$, then $\Lambda_{g_t}$ is a horseshoe. Thus, it remains only to estimate the Lebesgue measure of the set of parameters $t$ such that $d(A_{ct},B_{ct})<2ct$ in order to deduce that $\Lambda_{g}$ is a horseshoe for most $g\in\mathcal{U}$. At this stage, one notes that $A_{ct}$ is the $c\cdot t$-neighborhood of $K^s(g_t)=K^s(g_0)$ and $B_{ct}$ is the $c\cdot t$-neighborhood of $K^u(g_t)=K^s(g_0)+t$, and one recalls that, by Proposition 3 at page 104 of \cite{PT}, the fact that the regular Cantor sets $K^s(g_0)$ and $K^u(g_0)$ satisfy $d_s^0+d_u^0<1$ implies that for all $\varepsilon>0$ there exists $c=c(\varepsilon)>0$ and $t(\varepsilon)>0$ such that 
$$\frac{\textrm{Leb}(\{s\in(0,t): d(A_{cs},B_{cs})<2cs\})}{t}<\varepsilon$$
for each $0<t<t(\varepsilon)$. Of course, the reader has no difficulty to recognize that this last estimate readily implies that $\Lambda_g$ is horseshoe for most $g\in\mathcal{U}_0$, and thus the sketch of proof of Theorem \ref{t.Newhouse-Palis-Takens} is complete. 

After this discussion of homoclinic bifurcations of quadratic tangencies associated to thin horseshoes, we'll dedicate the rest of these notes to the study of bifurcations associated to fat horseshoes. 

\subsection{Homoclinic bifurcations associated to fat horseshoes and stable tangencies} Consider again the setting (and notations) of Subsection \ref{ss.bifurcations-Newhouse}: $f:M\to M$ is a $C^2$-diffeomorphism of a surface $M$ with a horseshoe $K$ and a periodic point $p\in K$ whose stable and unstable manifolds have a quadratic tangency at some point $q\in M-K$. Consider again $U$ a sufficiently small neighborhood of $K$ and $V$ a sufficiently small neighborhood of the orbit of $q$, and let's fix $\mathcal{U}$ a sufficiently small $C^2$-neighborhood of $f$ organized into the open sets $\mathcal{U}_-$ and $\mathcal{U}_+$ and the codimension $1$ hypersurface $\mathcal{U}_0$ depending on the relative positions of $W^s(p_g)$ and $W^u(p_g)$ near $V$. 

Assume further that the quadratic tangency of $f=g_0\in\mathcal{U}_0$ is associated to a \emph{fat} horseshoe\footnote{Here, we are implicitly using the continuity of the Hausdorff dimension of horseshoes to ensure that, if we choose $\mathcal{U}$ sufficiently small, then $\textrm{HD}(K_g)>1$ for $g\in\mathcal{U}$ once we have $\textrm{HD}(K)>1$ for some $f\in\mathcal{U}_0$. See \cite{PV} for more details on this.} $K=K_f$, i.e., suppose that $\textrm{HD}(K)=d_s^0+d_u^0>1$. By reviewing the arguments in the previous subsection, we see that the intersections or arithmetic differences of the regular Cantor sets $K^s(g_0)$ and $K^u(g_0)$ will hint what we should expect for the local dynamics of $g\in\mathcal{U}_+$. 

Here, the following result due to J. M. Marstrand is very inspiring: 

\begin{proposition}[J. M. Marstrand]\label{p.Marstrand} Let $C\subset\mathbb{R}^2$ be a subset with $\textrm{HD}(C)>1$. Then, for Lebesgue almost every $\lambda\in\mathbb{R}$, the set $\pi_{\lambda}(C)\subset \mathbb{R}$ has positive Lebesgue measure, where $\pi_{\lambda}(x,y):=x-\lambda y$.
\end{proposition}

For a proof of this result using potential theory, see Theorem 2 at page 64 of \cite{PT}.

In our context, we can apply Marstrand's theorem to $C=K^s(g_0)\times K^u(g_0)$ because, by hypothesis, $\textrm{HD}(C)\geq d_s^0+d_u^0>1$. By doing so, we get that for Lebesgue almost every $\lambda\in\mathbb{R}$ the arithmetic difference $K^s(g_0)\ominus\lambda K^u(g_0) = \pi_{\lambda}(K^s(g_0)\times K^u(g_0))$ has positive Lebesgue measure. In particular, if one can produce a $2$-parameter family $g_{\lambda,t}\in\mathcal{U}$ such that $g_{\lambda,0}\in\mathcal{U}_0$, $K^s(g_{\lambda,t})=K^s(g_0)$ and $K^u(g_{\lambda,t})=\lambda K^u(g_0) + t$ for all $\lambda$ close to $1$, then one would get that for almost every $\lambda$ close to $1$, the stable and unstable laminations of $K_{g_{\lambda,t}}$ meet tangentially in the region $V$ near $q$ for a set of parameters $t$ of positive Lebesgue measure. 

In particular, since the presence of tangencies prevents hyperbolicity, this hints that, in the context of fat horseshoes $K=K_f$, $f\in\mathcal{U}_0$, the statement of the theorem of Newhouse, Palis and Takens that $\Lambda_g$ is a horseshoe for most $g\in\mathcal{U}_+$ (cf. Theorem \ref{t.Newhouse-Palis-Takens}) may fail along certain $2$-parameter families $g_{\lambda,t}\in\mathcal{U}$. 

This idea was pursued in the work \cite{PY94} where J. Palis and J.-C. Yoccoz showed (in 1994) the following result. Let 
$$\mathcal{T}:=\{g\in\mathcal{U}: \textrm{the stable and unstable laminations of } K_g \textrm{ meet tangentially somewhere in } V\}$$
be the locus of tangencies. Then, for any smooth $2$-parameter family 
$(g_{\lambda,t})_{|\lambda|<\lambda_0, |t|<t_0}\in\mathcal{U}$ such that, for each $|\lambda|<\lambda_0$, 
\begin{itemize}
\item $f_{\lambda}=g_{\lambda,0}\in\mathcal{U}_0$, $g_{\lambda,t}\in\mathcal{U}_-$ for all $-t_0<t<0$ and $g_{\lambda,t}\in\mathcal{U_+}$ for all $0<t<t_0$, 
\item $g_{\lambda,t}$ is transverse to $\mathcal{U}_0$ at $f_{\lambda}:=g_{\lambda,0}$, 
\item the eigenvalues $\alpha_\lambda$ and $\beta_\lambda$ of the derivative of $f_{\lambda}:=g_{\lambda,0}$ at the periodic point $p_{\lambda,0}=p_{g_{\lambda,0}}$ vary \emph{non-trivially} with the parameter $\lambda$ in the sense that the derivatives $\frac{d\alpha_{\lambda}}{d\lambda}(0)$ and $\frac{d\beta_{\lambda}}{d\lambda}(0)$ at $\lambda=0$ do not vanish,
\item the horseshoe $K_{f_\lambda}$ is fat (i.e., its Hausdorff dimension is larger than 1),
\end{itemize}
there exists a constant $c=c(g_{\lambda,t})>0$ with the following property. For Lebesgue almost every $|\lambda|<\lambda_0$, it holds 
$$\limsup\limits_{\varepsilon\to0}\frac{\textrm{Leb}(\{t\in(0,t_0): g_{\lambda,t}\in\mathcal{T}\})}{\varepsilon}>c.$$

More recently, C. G. Moreira and J.-C. Yoccoz studied in \cite{MY01} the geometry of the intersections $K\cap K'$ of regular Cantor sets $K$ and $K'$, and they showed in \cite{MY10} how the key ideas from \cite{MY01} can be extended (with some non-trivial technical work) to show that the subset of $g\in\mathcal{U}_+$ with stable tangencies in the region $V$ (near $q$) has positive density in the setting of bifurcations of fat horseshoes. More precisely, let $\textrm{int}(\mathcal{T})$ be the locus of \emph{stable tangencies}, that is, $\textrm{int}(\mathcal{T})$ is the interior of $\mathcal{T}$. 

\begin{theorem}[C. G. Moreira and J.-C. Yoccoz]\label{t.MY10} Suppose that $\textrm{HD}(K)>1$ for $f\in\mathcal{U}_0$. Then, there exists an open and dense subset 
$\mathcal{U}_0^*$ of $\mathcal{U}_0$ such that any smooth $1$-parameter family $(g_t)_{|t|<t_0}$ passing by $g_0\in\mathcal{U}_0^*$ transversely to $\mathcal{U}_0$ meets the locus of stable tangencies $\textrm{int}(\mathcal{T})$ with positive (inferior) density in the sense that 
$$\liminf\limits_{\varepsilon\to0}\frac{\textrm{Leb}(\{t\in(0,\varepsilon):\, g_t\in\textrm{int}(\mathcal{T})\})}{\varepsilon}>0$$
Furthermore, denoting by $\mathcal{H}:=\{g\in\mathcal{U}:\,\Lambda_g \textrm{ is a horseshoe}\}$, one has that $(g_t)_{|t|<t_0}$ meets $\textrm{int}(\mathcal{T})\cup\mathcal{H}$ with full density in the sense that 
$$\lim\limits_{\varepsilon\to0}\frac{\textrm{Leb}(\{t\in(0,\varepsilon):\, g_t\in\textrm{int}(\mathcal{T})\cup\mathcal{H}\})}{\varepsilon}=1$$
\end{theorem}

Evidently, this result makes clear that in the context of fat horseshoes (i.e., $\textrm{HD}(K)>1$ for $f\in\mathcal{U}_0$) it is not true that $\Lambda_g$ is a horseshoe for most $g\in\mathcal{U}_+$. 

We can summarize the discussion so far with the following two phrases:
\begin{itemize}
\item by Theorem \ref{t.Newhouse-Palis-Takens}, in the context of thin horseshoes (i.e., $\textrm{HD}(K)<1$ for $f\in\mathcal{U}_0$), $\Lambda_g$ is a horseshoe for most $g\in\mathcal{U}_+$, and 
\item by Theorem \ref{t.MY10}, in the context of fat horseshoes (i.e., $\textrm{HD}(K)>1$ for $f\in\mathcal{U}_0$), $\Lambda_g$ has persistent tangencies with positive ``probability'' and thus we can't expect that $\Lambda_g$ is a horseshoe for most $g\in\mathcal{U}_+$.
\end{itemize}

\subsection{Heteroclinic bifurcations of slightly fat horseshoes after J. Palis and J.-C. Yoccoz} We saw that bifurcations of quadratic tangencies associated to fat horseshoes, 
$\textrm{HD}(K)=d_s^0+d_u^0>1$, are complicated because of persistent tangencies. However, by a closer inspection of the works \cite{MY01} and \cite{MY10}, one realizes that the regular Cantor sets $K^s(g)$ and $K^u(g)$ for $g\in\mathcal{U}_+$ are usually expected to intersect in a set $K^s(g)\cap K^u(g)$ of Hausdorff dimension \emph{close} to $d_s^0+d_u^0-1=\textrm{HD}(K)-1$. Thus, from the heuristic point of view, the \emph{critical locus} $K^s(g)\cap K^u(g)$ (i.e., the region where the tangencies destroying the hyperbolicity show up) is very \emph{small}, i.e., its Hausdorff dimension is close to zero, if the initial horseshoe $K$ is only \emph{slightly} fat, i.e., $\textrm{HD}(K)>1$ is close to $1$. In particular, one could imagine that bifurcations quadratic tangencies of slightly fat horseshoes could lead to a local dynamics on $\Lambda_g$ satisfying some form of weak (non-uniform) hyperbolicity for most $g\in\mathcal{U}_+$ despite the fact that $\Lambda_g$ doesn't verify strong (uniform) hyperbolicity conditions in general. 

In a recent tour-de-force work (of 217 pages), J. Palis and J.-C. Yoccoz \cite{PY} were able to formalize this crude heuristic argument by showing (among several other things) the following result in the context of heteroclinic bifurcations of slightly fat horseshoes. 

Let $f$ be a smooth diffeomorphism of a compact surface $M$ possessing a uniformly hyperbolic horseshoe $K$ displaying a \emph{heteroclinic} quadratic tangency, that is,  $K$ contains two periodic points $p_s$ and $p_u$ with distinct orbits such that $W^s(p_s)$ and $W^u(p_u)$ have a quadratic tangency (i.e., a contact of order $1$) at some point $q\in M-K$. Let $U$ be a sufficiently small neighborhood of $K$ and let $V$ be a sufficiently small neighborhood of $q$ such that $K\cup\mathcal{O}(q)$ is the maximal invariant set of $U\cup V$. Denote by $\mathcal{U}$ be a sufficiently small neighborhood of $f$ and, as usual, let's organize $\mathcal{U}$ into $\mathcal{U}=\mathcal{U}_-\cup\mathcal{U}_0\cup\mathcal{U}_+$ depending on the relative positions of the continuations of $W^s(p_s)$ and $W^u(p_u)$ near $V$ (see Figure \ref{f.4}). Finally, let us denote by $d_s^0$ and $d_u^0$ the stable and unstable dimensions of the horseshoe $K$ of $f\in\mathcal{U}_0$.

\begin{theorem}[J. Palis and J.-C. Yoccoz]\label{t.PY} In the setting of the paragraph above, suppose that $K$ is slightly fat in the sense that 
\begin{equation}\label{e.PYcondition}
(d_s^0+d_u^0)^2+(\max\{d_s^0,d_u^0\})^2 < (d_s^0+d_u^0)+(\max\{d_s^0,d_u^0\})
\end{equation}
Then, $\Lambda_g$ is a non-uniformly hyperbolic horseshoe for most $g\in\mathcal{U}_+$.
\end{theorem} 

\begin{remark} At first sight, there is no reason to restrict our attention to heteroclinic tangencies in the previous theorem. In fact, as we'll see later (cf. Remark \ref{r.PYheteroclinic}), for certain technical reasons, the arguments of J. Palis and J.-C. Yoccoz can treat only heteroclinic tangencies. Of course, the authors believe that this is merely an artifact of their methods, but unfortunately they don't know how to modify the proofs to also include the case of homoclinic tangencies.
\end{remark}

Concerning the statement of this result, let us comment first on condition \eqref{e.PYcondition}. As a trivial remark, note that this condition includes the case $d_s^0+d_u^0<1$ of thin horseshoes, but this is not surprising as any reasonable definition of ``non-uniformly hyperbolic horseshoe'' must include uniformly hyperbolic horseshoes as particular examples. Of course, this remark is not particularly interesting because the case of thin horseshoes was already treated by S. Newhouse, J. Palis and F. Takens (cf. Theorem \ref{t.Newhouse-Palis-Takens}), so that condition \eqref{e.PYcondition} is really interesting in the regime of fat horseshoes $d_s^0+d_u^0>1$. Here, one can get a clear idea about \eqref{e.PYcondition} by assuming $\max\{d_s^0,d_u^0\}=d_s^0$ or $d_u^0$ (i.e., by breaking the natural symmetry between $d_s^0$ and $d_u^0$), and by noticing that the boundary of the region determined by \eqref{e.PYcondition} is the union of two ellipses meeting the diagonal $\{d_s^0=d_u^0\}$ at the point $(3/5,3/5)$ as indicated in the figure below:

\begin{figure}[htb!]
\includegraphics[scale=0.48]{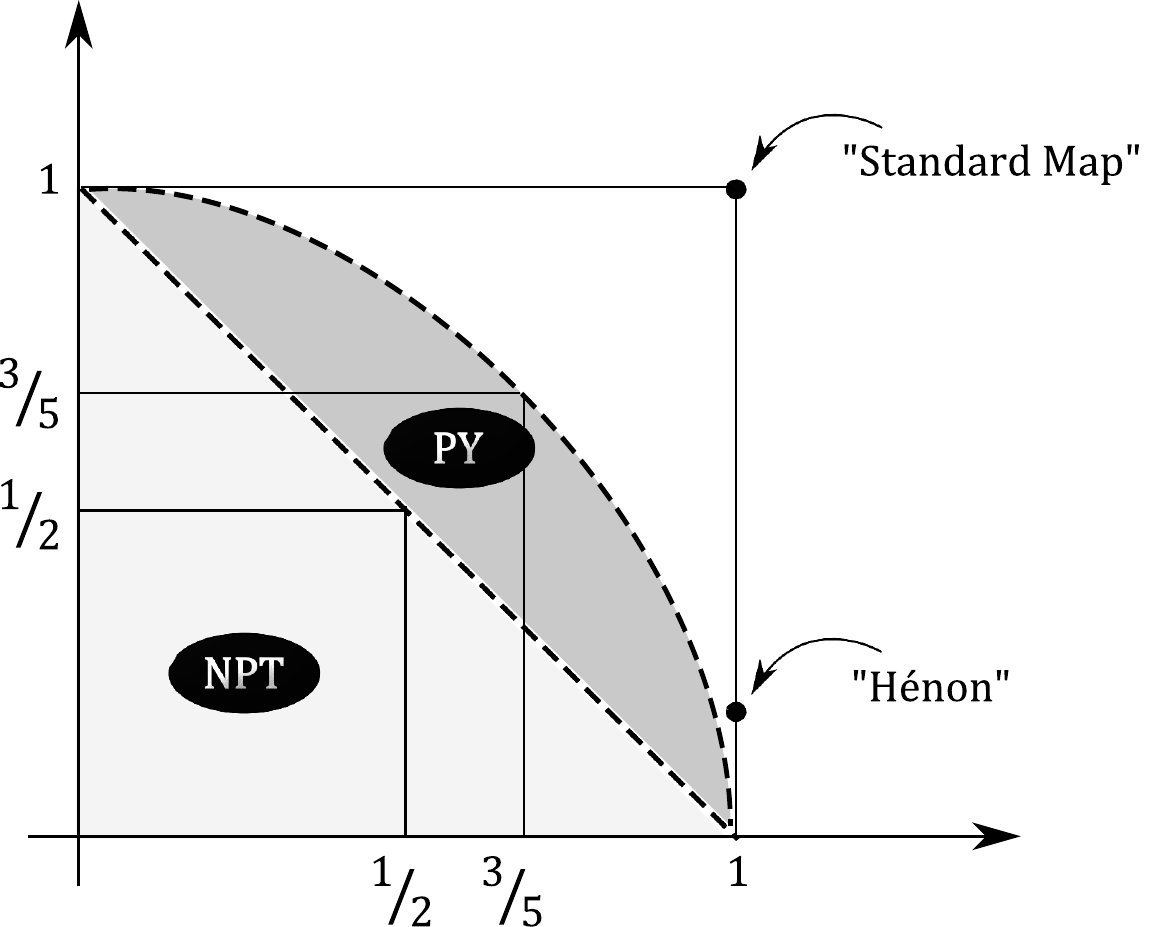}
\caption{Region of parameters $d_s^0$ and $d_u^0$ where the results of Newhouse-Palis-Takens (NPT) and Palis-Yoccoz (PY) apply.}\label{f.10}
\end{figure}

In this figure, we used the horizontal axis for the variable $d_s^0$ and the vertical axis for the variable $d_u^0$. Also, we pointed out, for sake of comparison, two famous families of dynamical systems lying outside the scope of Theorem \ref{t.PY}, namely the \emph{H\'enon maps} $H_{a,b}:\mathbb{R}^2\to\mathbb{R}^2$, $H_{a,b}(x,y)=(1-ax^2+y,bx)$, and the \emph{standard family} $f_{\lambda}:\mathbb{T}^2\to\mathbb{T}^2$, $f_{\lambda}(x,y)=(2x+\lambda\sin(2\pi x)-y, x)$. Indeed, these important examples of dynamical systems can not be studied by the current methods of J. Palis and J.-C. Yoccoz because they display homoclinic and heteroclinic bifurcations associated to ``very fat horseshoes'': 
\begin{itemize}
\item in the case of H\'enon maps, the horseshoes have stable dimension $d_s^0=1$ and a very small unstable dimension $0<d_u^0\ll1$ for certain parameters $(a,b)$, and 
\item in the case of the standard family, one has horseshoes with $d_s^0=d_u^0$ arbitrarily close to $1$ for large values of $\lambda\in\mathbb{R}$.
\end{itemize}

For further comments on other works on non-uniform hyperbolicity results inside parametrized families (such as H\'enon maps), we refer the reader to the book \cite{PT} of J. Palis and F. Takens (and the references therein).

Now, let us start to explain the meaning of \emph{non-uniformly hyperbolic horseshoe} in Theorem \ref{t.PY}. As we explained in Definition \ref{d.UH-horseshoe}, a (uniformly hyperbolic) horseshoe $\Lambda$ of a surface diffeomorphism $f:M\to M$ is a \emph{saddle-like} object in the sense that $\Lambda$ is not an attractor nor a repellor, that is, both its \emph{stable set}
$$W^s(\Lambda):=\{y\in M: \textrm{dist}(f^n(y),\Lambda)\to 0 \textrm{ as }n\to+\infty\}$$
and \emph{unstable set}
$$W^u(\Lambda):=\{y\in M: \textrm{dist}(f^n(y),\Lambda)\to 0 \textrm{ as }n\to+\infty\}$$
have zero Lebesgue measure $\textrm{Leb}_2$. Here, $\textrm{Leb}_2$ is the $2$-dimensional Lebesgue measure of $M$. In a similar vein, J. Palis and J.-C. Yoccoz (cf. Theorem 7 of \cite{PY}) showed that their non-uniformly hyperbolic horseshoes are saddle-like objects:

\begin{theorem}\label{t.PY-thm7} Under the same assumptions of Theorem \ref{t.PY}, one has that 
$$\textrm{Leb}_2(W^s(\Lambda_g))=\textrm{Leb}_2(W^u(\Lambda_g))=0$$ 
for most $g\in\mathcal{U}_+$.
\end{theorem}

A nice way to better appreciate this statement is to contrast it with Newhouse phenomena (cf. Theorem~\ref{t.Newhouse} and Remark~\ref{r.Newhouse-w/o-thickness}). Indeed, while Newhouse phenomena ensure that the coexistence of infinitely many sinks/sources inside $\Lambda_g$ for some $g\in\mathcal{U}_+$, we know from Theorem~\ref{t.PY-thm7} that $\Lambda_g$ doesn't contains sinks or sources for most $g\in\mathcal{U}_+$.

Actually, the statement of Theorem 7 of \cite{PY} contains a slightly more precise explanation of the non-uniformly hyperbolic features of $\Lambda_g$ (for most 
$g\in\mathcal{U}_+$): it is possible to show that $\Lambda_g$ supports geometric \emph{Sinai-Ruelle-Bowen} (SRB) measures with \emph{non-zero Lyapunov exponents}, that is, $\Lambda_g$ is a non-uniformly hyperbolic object in the sense of the so-called \emph{Pesin theory}. Unfortunately, a detailed explanation of these terms (i.e., SRB measures, Lyapunov exponents, Pesin theory) is out of the scope of these notes and we refer the curious reader to \cite{BR}, \cite{Si} and \cite{KH} for more informations.

In order to further explain the structure of $\Lambda_g$, we will briefly describe in the next subsection some elements of the proof of Theorem \ref{t.PY}.

\subsection{A global view on Palis-Yoccoz induction scheme}\label{ss.PY-global-view} Let $(g_t)_{|t|<t_0}$ a smooth $1$-parameter family transverse to $\mathcal{U}_0$ at $f=g_0$, where $f$ is a diffeomorphism with a slightly fat horseshoe $K$ exhibiting a heteroclinic quadratic tangency as shown in the figure below:
\begin{figure}[htb!]
\includegraphics[scale=0.47]{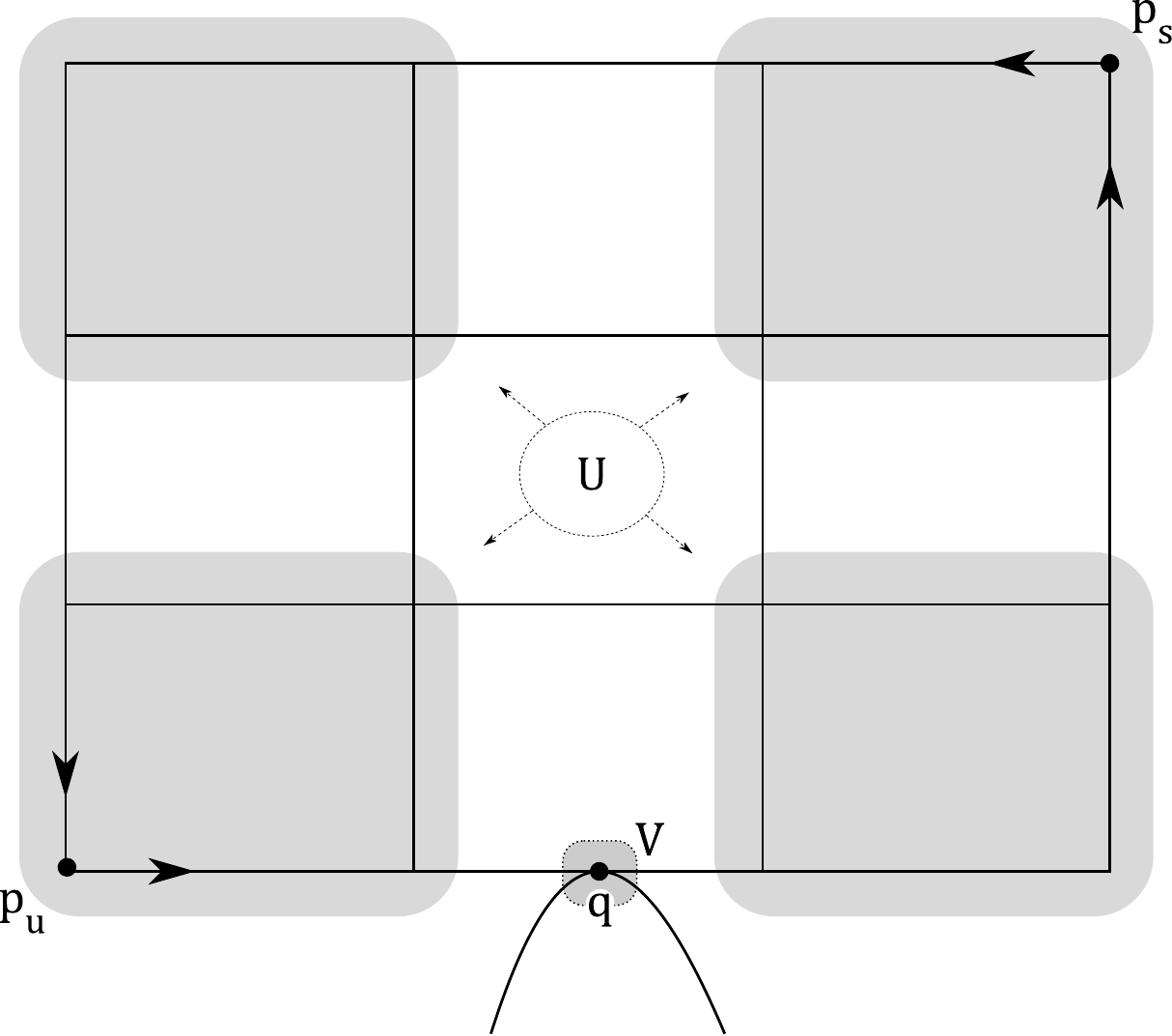}
\caption{Heteroclinic quadratic tangency associated to a slightly fat horseshoe.}\label{f.11}
\end{figure}

As usual, we wish to understand the local dynamics of $g=g_t$ on the neighborhood $U\cup V$ indicated in the picture above, that is, we want to investigate the structure of the set  
$$\Lambda_g=\bigcap\limits_{n\in\mathbb{Z}} g^n(U\cup V)$$
for most parameters $t\in (0,t_0)$. 

In this direction, we consider $0<\varepsilon_0\ll1$ and we look at the parameter interval $I_0:=[\varepsilon_0,2\varepsilon_0]$. Very roughly speaking, the scheme of J. Palis and J.-C. Yoccoz has the following structure: besides $\varepsilon_0$, one has two extra parameters $\tau$ and $\eta$ chosen such that $$0<\varepsilon_0\ll\eta\ll\tau\ll1$$
Then, one proceeds inductively: 
\begin{itemize}
\item at the 1st stage, one defines $\varepsilon_1:=\varepsilon_0^{1+\tau}$ and one divides the interval $I_0:=[\varepsilon_0, 2\varepsilon_0]$ into $[\varepsilon_0^{-\tau}]$ candidate subintervals;
\item then, one apply an \emph{exam} called \emph{strong regularity} to each candidate subinterval: the \emph{good} subintervals (those passing the strong regularity test) are \emph{kept} while the bad ones are discarded; 
\item after that, one goes to the next stages, that is, one takes the good intervals $I_k$ from the $k$th stage, subdivide them into $[\varepsilon_k^{-\tau}]$ subintervals of size 
$\varepsilon_{k+1}:=\varepsilon_k^{1+\tau}$, apply the strong regularity exam to each subinterval and one keeps the good subintervals and discard the bad subintervals.
\end{itemize}

Of course, the strong regularity of an interval $I$ is a property about the (non-uniform) hyperbolic features of $\Lambda_{g_t}$ for all parameters $t\in I$, and the choice of the set of properties defining the strong regularity must be extremely careful: it should not be too weak (otherwise one doesn't get hyperbolicity) nor too strong (otherwise there is a risk that no interval is good at some stage).

Actually, as we'll see later, for each candidate interval $I$, J. Palis and J.-C. Yoccoz construct a class $\mathcal{R}(I)$ of so-called ($I$-persistent) affine-like iterates\footnote{Roughly speaking, these are iterates of $g$ behaving like affine maps of the plane in the sense of Definitions~\ref{d.weak-affine-like} and~\ref{d.affine-like} below.} of $g=g_t$, $t\in I$ and they test the strong regularity of $I$ by examining the features of the class $\mathcal{R}(I)$. 

\begin{remark}\label{r.past-future} A nice feature of the arguments of J. Palis and J.-C. Yoccoz is that they are \emph{time-symmetric}, that is, the dynamical estimates for the past and the future are the \emph{same} (i.e., one has only to do half of the computations). In particular, those readers with some familiarity with H\'enon maps know that the past behavior is very different from the future behavior (due to strong dissipation) and this essentially explains why the methods of \cite{PY} are not directly useful in the case of H\'enon maps.
\end{remark}

After this very approximative description of Palis-Yoccoz inductive scheme, it is clear that one of the key ideas is to carefully setup the notion of strong regularity property. However, before discussing this subject, we need to make some preparations: firstly we need to localize the dynamics, secondly we need to introduce the affine-like iterates, and thirdly we need to introduce the class $\mathcal{R}(I)$.

\subsubsection{Localization of the dynamics}

The local dynamics of $g_t$ for $t\in I_0:=[\varepsilon_0,2\varepsilon_0]$ has the following appearance:
\begin{figure}[htb!]
\includegraphics[scale=0.7]{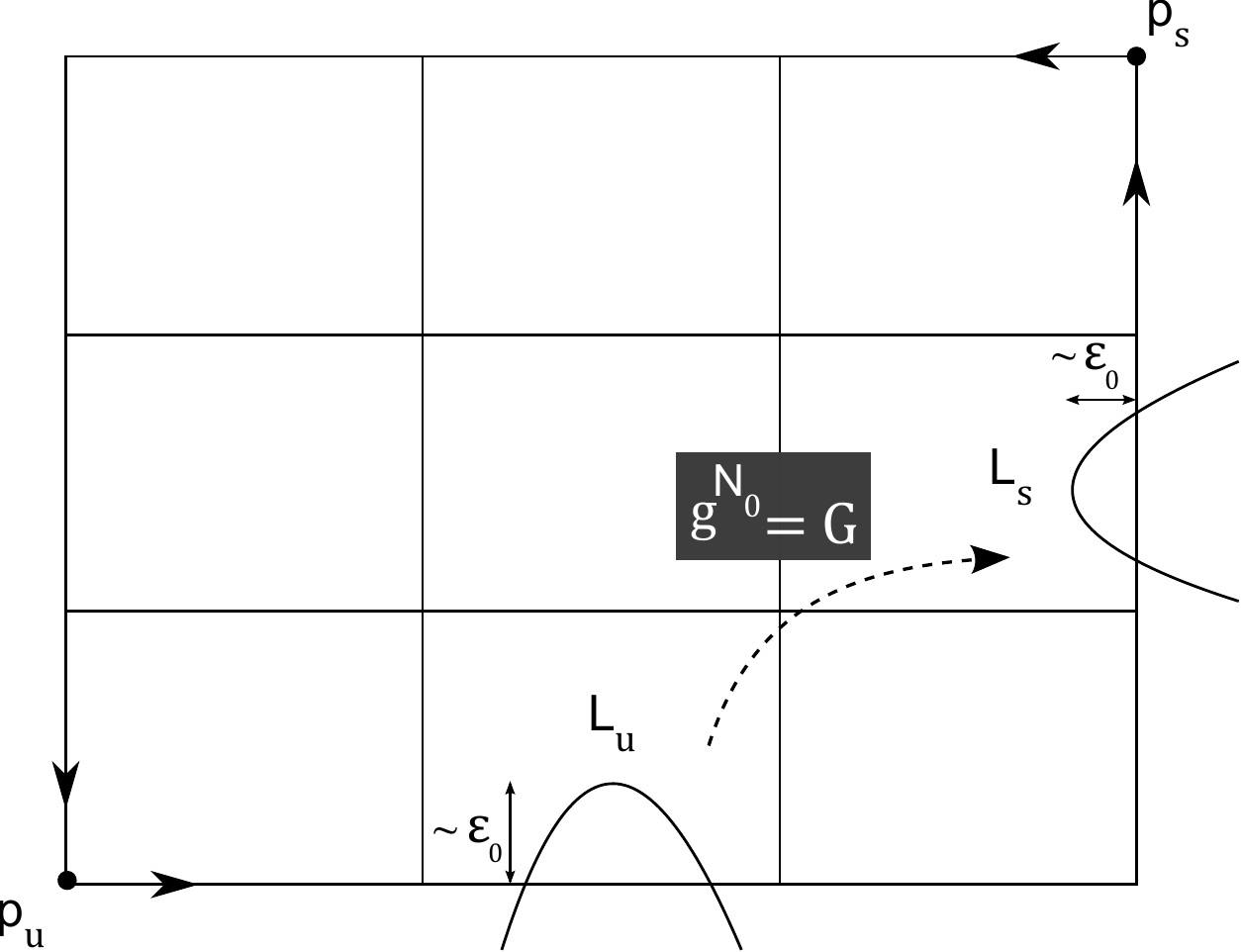}
\caption{Parabolic tongues created after unfolding a heteroclinic tangency.}\label{f.12}
\end{figure}

As it is highlighted in this picture, after unfolding the tangency, we get two regions $L_u$ and $L_s$ called \emph{unstable} and \emph{stable} \emph{parabolic tongues} bounded by the pieces of $W_{loc}^s(p_s)$ and $W^u_{loc}(p_u)$ near $V$. The transition time $N_0$ from the unstable tongue $L_u$ to the stable tongue $L_s$ under the dynamics $g$ is a large but fixed integer depending only on $f\in\mathcal{U}_0$. 

Using this, we can organize the local dynamics of $g$ on $U\cup V$ as follows. Firstly, we select a finite Markov partition of the horseshoe $K_g=\bigcap\limits_{n\in\mathbb{Z}}g^n(U)$ into compact disjoint rectangles $R_a$, $a\in\mathcal{A}$, i.e., by fixing a convenient system of coordinates, we write $R_a\simeq I_a^u\times I_a^s$ in such a way that: 
\begin{itemize}
\item the derivative of $g|{R_a}$ expands (uniformly) the horizontal direction and contracts (uniformly) the vertical direction, 
\item $K_g$ is the maximal invariant set of the interior $\textrm{int}(R)$ of $R:=\bigcup\limits_{a\in\mathcal{A}} R_a$, i.e., $K_g=\bigcap\limits_{n\in\mathbb{Z}}g^n(\textrm{int}(R))$, 
\item $g(I_a^u\times\partial I_a^s)\cap \textrm{int}(R)=\emptyset=g^{-1}(\partial I_a^u\times I_a^s)\cap \textrm{int}(R)$;
\item $\{R_a\cap K_g\}_{a\in\mathcal{A}}$ is a \emph{Markov partition}, i.e., $g^{-1}(I_a^u\times\partial I_a^s)\subset \bigcup\limits_{b\in\mathcal{A}}(I_b^u\times\partial I_b^s)$, $g(\partial I_a^u\times I_a^s)\subset \bigcup\limits_{b\in\mathcal{A}}(\partial I_b^u\times I_b^s)$, and there exists an integer $n\in\mathbb{N}$ with $g^n(R_a)\cap R_{a'}\neq\emptyset$ for all $a,a'\in\mathcal{A}$.
\end{itemize}
Secondly, we denote by $\mathcal{B}=\{(a,a')\in\mathcal{A}^2: g(R_a)\cap R_{a'}\neq\emptyset\}$. Then, in this setting, it is not hard to see that the local dynamics of $g$ on $U\cup V$ is given by 
\begin{itemize}
\item the \emph{uniformly hyperbolic maps}\footnote{They are called uniformly hyperbolic because the horizontal direction is uniformly expanded and the vertical direction is uniformly contracted by their derivatives.} $g:R_a\cap g^{-1}(R_a)\to g(R_a)\cap R_{a'}$, $(a,a')\in\mathcal{B}$ related to the horseshoe $K_g$, and
\item the \emph{folding map} $G:=g^{N_0}: L_u\to L_s$ making the transition between parabolic tongues.
\end{itemize}

In this context, by letting $\widehat{R}=R\cup\bigcup\limits_{i=1}^{N_0-1}g^i(L_u)$, we have that $\Lambda_g$ is the maximal invariant set of $\widehat{R}$, i.e., $\Lambda_g=\bigcap\limits_{n\in\mathbb{Z}}g^n(\widehat{R})$.

This localization of the dynamics of $g$ on $\Lambda_g$ to the region $\widehat{R}$ is useful because it allows us to think of $\Lambda_g$ in terms of an \emph{iterated system of maps}, i.e., we approach the points of $\Lambda_g$ by looking at the domains and the images of the \emph{compositions} (i.e., certain $g$-iterates) of the uniformly hyperbolic maps $g:R_a\cap g^{-1}(R_a)\to g(R_a)\cap R_{a'}$ and the folding map $G=g^{N_0}: L_u\to L_s$. 

By thinking in this way, we see that the points in the domains or images of $g$-iterates (composition) with \emph{affine-like features}, that is, $g$-iterates whose derivates expand the horizontal direction and contract the vertical direction, will contribute to the hyperbolicity of $\Lambda_g$. In other words, it is desirable to get as much affine-like iterates as possible. Of course, the $g$-iterates obtained by composition of transition maps $g:R_a\cap g^{-1}(R_a)\to g(R_a)\cap R_{a'}$ related to the horseshoe $K_g$ have affine-like features (by definition), so that one risks losing the affine-like property only when one considers compositions with the folding map $G$ (because the folding map mixes up the horizontal and vertical directions). 

In particular, this suggests that the strong regularity property has something to do with the consecutive passages through the critical region given by the parabolic tongues $L_u$ and $L_s$. However, before pursuing this direction, let us formalize the notion of \emph{affine-like iterates}.

\subsubsection{Affine-like maps} A \emph{vertical strip} $P\subset R_0\simeq I_0^u\times I_0^s$ is a region of the form 
$$P=\{(x_0,y_0)\in R_0: \phi^-(y_0)\leq x_0\leq \phi^+(y_0)\}$$
and a \emph{horizontal strip} $Q\subset R_1\simeq I_1^u\times I_1^s$ is a region of the form
$$Q=\{(x_1,y_1)\in R_1: \psi^-(x_1)\leq y_1\leq \psi^+(x_1)\}$$ 
Intuitively, we wish to call ``affine-like'' a map $F:P\to Q$ from a vertical strip $P$ to a horizontal strip $Q$ approximately contracting the vertical direction and expanding horizontal direction such as the one depicted in Figure \ref{f.13} below.
\begin{figure}[htb!]
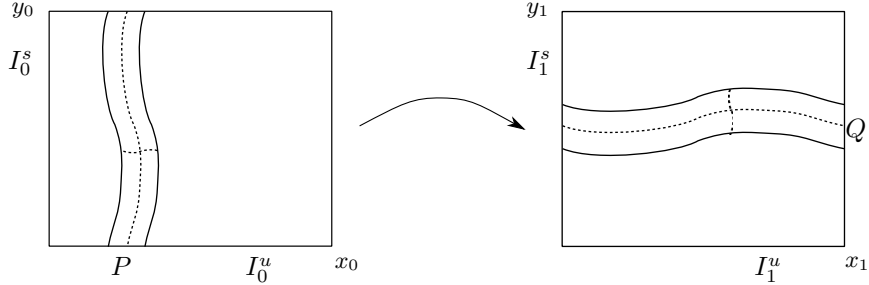
\caption{Geometry of affine-like maps.}\label{f.13}
\end{figure}

Formally, we define:
\begin{definition}\label{d.weak-affine-like}We say that a map $F(x_0,y_0)=(x_1,y_1)$ from a vertical strip $P$ to a horizontal strip $Q$ is \emph{weakly affine-like} whenever $F$ admits an \emph{implicit representation} $(A,B)$, i.e., we can write $x_0=A(x_1,y_0)$ and $y_1=B(x_1,y_0)$. Equivalently, $F$ is weakly affine-like if and only if the projection $\pi$ from the graph $\textrm{graph}(F)$ of $F$ to $I_0^u\times I_1^s$ is a diffeomorphism.
\end{definition}

This definition of affine-like maps $F$ in terms of implicit representations $(A,B)$ is somewhat folkloric in Dynamical Systems, and it was used by J. Palis and J.-C. Yoccoz because it is technically easier to estimate $(A,B)$ than $F$ as the symmetry between past and future is more evident, and $(A,B)$ are contractive maps. 

In what follows, we will denote the derivatives of $A$ and $B$ by $A_x, B_x, A_y, B_y, A_{xx}, B_{xx}, \dots$. Following J. Palis and J.-C. Yoccoz, we will consider \emph{exclusively} weakly affine-like maps satisfying the following \emph{hyperbolicity} conditions:
\begin{definition}\label{d.affine-like} A weakly affine-like map $F$ is called \emph{affine-like} if its implicit representation $(A,B)$ verifies:
\begin{itemize}
\item \emph{Cone condition}: $\lambda|A_x|+u|A_y|\leq 1$ and $\lambda|B_y|+v|B_x|\leq 1$ where $1<uv<\lambda^2$, and
\item \emph{Bounded distortion condition}: $\partial_x\log |A_x|, \partial_y\log |A_x|, A_{yy}, \partial_y\log |B_y|, \partial_x\log |B_y|, B_{xx}$ are uniformly bounded by some constant $C>0$.
\end{itemize}
Here, the constants $\lambda, u, v$ and $C$ are fixed once and for all depending only on $f=g_0\in\mathcal{U}_0$.
\end{definition}

Informally, the cone condition says that $F$ contracts the vertical direction and expands the horizontal direction, and the bounded distortion condition says that the derivative of 
$F$ behaves in the same way in all scales.

For later use, we introduce the following notion:
\begin{definition}The widths of the domain $P$ and the image $Q$ of an affine-like map $F:P\to Q$ with implicit representation $(A,B)$ are 
$$|P|=\max|A_x| \quad \textrm{ and } \quad |Q|=\max|B_y|$$
\end{definition}

Once we dispose of the notion of affine-like iterates, we're ready to introduce the class $\mathcal{R}(I)$ whose strong regularity will be tested later. 

\subsubsection{Simple and parabolic compositions of affine-like maps and the class $\mathcal{R}(I)$} Coming back to the interpretation of the dynamics on $\Lambda_g$ as an iterated system of maps given by compositions of $g:R_a\cap g^{-1}(R_a)\to g(R_a)\cap R_{a'}$ and the folding map $G:L_u\to L_s$, we see that the following two ways of composing affine-like maps are particularly interesting in our context. 

\begin{definition}Let $F:P\to Q$ and $F':P'\to Q'$ be two affine-like maps such that $Q, P'\subset R_{a'}$. Then, the \emph{simple composition} $F''=F'\circ F$ is the affine-like map with domain $P'':=P\cap F^{-1}(P')$ and image $Q'':=Q'\cap F'(Q)$ shown in Figure \ref{f.14} below.
\end{definition}
\begin{figure}[htb!]
\includegraphics[scale=0.63]{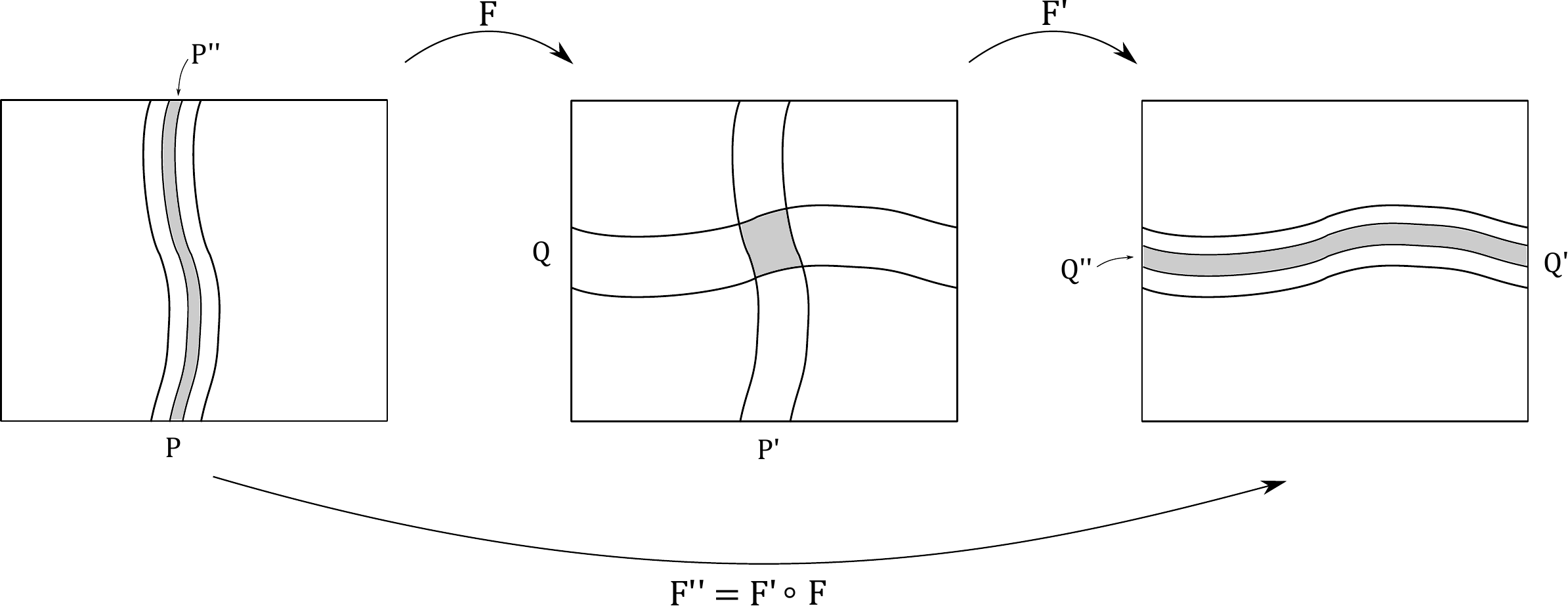}
\caption{Simple composition of affine-like maps.}\label{f.14}
\end{figure}
\begin{remark}
By direct inspection of definitions, one can check that $|P''|\sim|P|\cdot|P'|$ where the implied constant depends only on $f$ (by means of the constants $u,v$ in the cone condition).
\end{remark}

The composition of two transition maps $g:R_a\cap g^{-1}(R_a)\to g(R_a)\cap R_{a'}$ and $g:R_{a'}\cap g^{-1}(R_{a''})\to g(R_{a'})\cap R_{a''}$ associated to the horseshoe $K_g$ is the canonical example of simple composition. 

In particular, if we wish to understand $\Lambda_g$, it is not a good idea to work only with simple compositions, that is, we must include some passages through the parabolic tongues. This is formalized by the following notion. 

\begin{definition} Denote by $R_{a_u}$ and $R_{a_s}$ the rectangles of the Markov partition of $K_g$ containing the parabolic tongues $L_u$ and $L_s$. Let $F_0:P_0\to Q_0$ and $F_1=P_1\to Q_1$ be two affine-like maps such that $Q_0$, resp. $P_1$, passes near the parabolic tongue $L_u$, resp. $L_s$, i.e., $Q_0\subset R_{a_u}$ crosses $L_u$ and $P_1\subset R_{a_s}$ crosses $L_s$. We define the \emph{parabolic compositions} of $F_0$ and $F_1$ as follows. Firstly, we compare $Q_0$ with the parabolic-like strip $G^{-1}(P_1\cap L_s)$ and we say that the parabolic composition of $F_0$ and $F_1$ is \emph{possible} if the intersection $Q_0\cap G^{-1}(P_1\cap L_s)$ has two connected components $Q_0^-$ and $Q_0^+$ as shown (in black) in Figure \ref{f.15} below. Then, assuming that the parabolic composition of $F_0$ and $F_1$ is possible, we define their parabolic compositions as the two weakly affine-like maps $F^-:P^-\to Q^-$ and $F^+:P^+\to Q^+$ shown in Figure \ref{f.15} below obtaining by concatenating $F_0$, the folding map $G$ and $F_1$ in the strips $P^-=F_0^{-1}(Q_0^-)$, $P^+=F_0^{-1}(Q_0^+)$, $Q^-=F_1(G(Q_0^-))$, $Q^+=F_1(G(Q_0^+))$.
\end{definition}

\begin{figure}[htb!]
\includegraphics[scale=0.65]{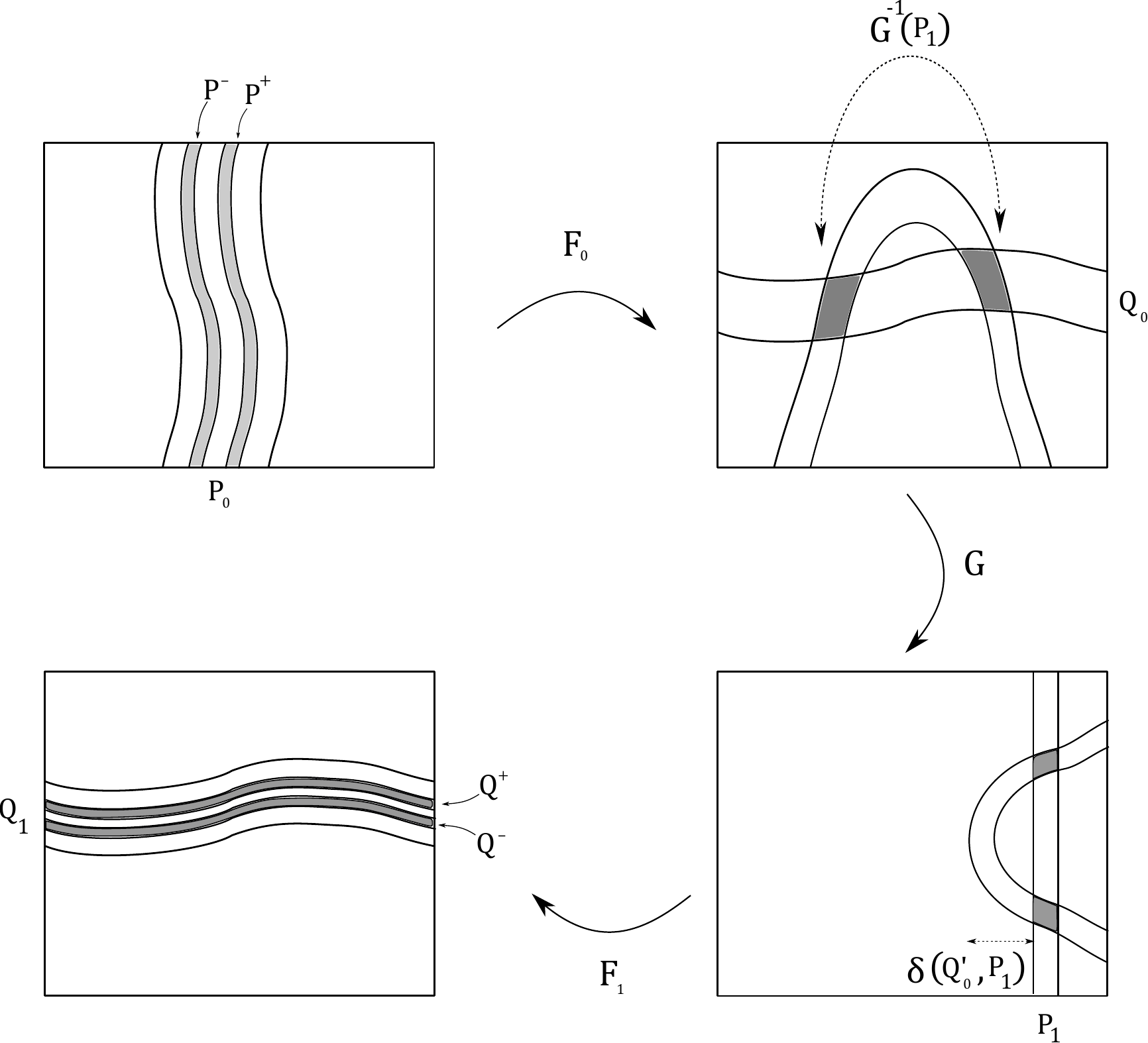}
\caption{Parabolic compositions of affine-like maps.}\label{f.15}
\end{figure}

As it is indicated in the figure above, the parabolic composition comes with an important parameter $\delta(Q_0,P_1)$ measuring the distance between the vertical strip $P_1$ and the tip of the parabolic-like strip $G(Q_0\cap L_u)$, or, equivalently, the horizontal strip $Q_0$ and the tip of the parabolic-like strip $G^{-1}(P_1\cap L_s)$.

\begin{remark} By direct inspection of definitions, one can check that $|P^{\pm}|=\frac{|P_0|\cdot|P_1|}{\delta(Q_0,P_1)^{1/2}}$.
\end{remark}

In this notation, the class $\mathcal{R}(I)$ is defined as follows.

\begin{definition}$\mathcal{R}(I)$ is the class of affine-like iterates of $g_t$, $t\in I$, closed under all simple compositions and certain parabolic compositions. More precisely, $\mathcal{R}(I)$ contains only parabolic compositions satisfying certain \emph{transversality} conditions such as 
$$\delta(Q_0,P_1)\geq \max\{|Q_0|^{1-\eta}, |P_1|^{1-\eta}, |I|\}.$$
\end{definition}

\begin{remark} In fact, the transversality conditions on parabolic compositions imposed by J. Palis and J.-C. Yoccoz involves 6 conditions besides the one on the parameter $\delta(Q_0, P_1)$ given above. Also, it is worth to point out that the class $\mathcal{R}(I)$ satisfying these conditions is unique, but this is shown in \cite{PY} only a posteriori.
\end{remark}

For later use, we denote by $(P,Q,n)$ an affine-like iterate $g^n:P\to Q$ taking a vertical strip $P$ to a horizontal strip $Q$ after $n$ iterations of $g=g_t$.

At this stage, we are ready to discuss the strong regularity property for $\mathcal{R}(I)$.

\subsubsection{Critical strips, bicritical dynamics and strong regularity} Let $(P,Q,n)\in\mathcal{R}(I)$.

\begin{definition}We say that $P$ is $I$-\emph{critical} when $P$ is not $I$-transverse to the parabolic tongue $L_s$, i.e., the distance between $P$ to the ``tip'' of $L_s$ is smaller $|P|^{1-\eta}$ for some $t\in I$. Similarly for $Q$ and $L_u$.
\end{definition}

\begin{definition}We say that an element $(P,Q,n)\in\mathcal{R}(I)$ is $I$-\emph{bicritical} if $P$ \emph{and} $Q$ are $I$-critical.
\end{definition}

In other words, a bicritical $(P,Q,n)\in\mathcal{R}(I)$ corresponds to some part of the dynamics starting at some $P$ close to the tip of $L_s$ and ending at some $Q$ close to the tip of $L_u$, that is, a bicritical $(P,Q,n)\in\mathcal{R}(I)$ corresponds to a return of the critical region to itself. 

Of course, one way of getting hyperbolicity for $\Lambda_g$ is to control the bicritical dynamics, i.e., bicritical elements $(P,Q,n)\in\mathcal{R}(I)$. 

\begin{definition} Given $\beta>1$, we say that a candidate parameter $I$ is $\beta$-\emph{regular} if 
$$|P|,|Q|<|I|^{\beta}$$
for every $I$-bicritical element $(P,Q,n)\in \mathcal{R}(I)$.
\end{definition}

\begin{remark}\label{r.beta} In their article, J. Palis and J.-C. Yoccoz choose $\beta>1$ depending only on the stable and unstable dimensions $d_s^0$ and $d_u^0$ of the initial horseshoe $K$ and the hyperbolicity strength of the periodic points $p_s$ and $p_u$ involved in the heteroclinic tangency. See Equation (5.19) of \cite{PY} for the precise requirements on 
$\beta$.
\end{remark}

Intuitively, a candidate parameter interval $I$ is $\beta$-regular if the bicritical dynamics seen through $\mathcal{R}(I)$ is confined to very small strips $P$ and $Q$. Unfortunately, the condition of $\beta$-regularity is not enough to run the induction scheme of J. Palis and J.-C. Yoccoz, and they end up by introducing a more technical condition called \emph{strong regularity}. However, for the sake of this text, we will assume that strong regularity \emph{is} $\beta$-regularity for some adequate parameter $\beta>1$. 

After this brief discussion of strong regularity, it is time to come back to Palis-Yoccoz induction scheme in order to say a few words about the dynamics of $\Lambda_{g_t}$ for $t$ belonging to strongly regular intervals.

\subsubsection{Dynamics of strong regular parameters} As it is explained in Sections 10 and 11 of \cite{PY}, J. Palis and J.-C. Yoccoz can reasonably control the dynamics of 
$\Lambda_{g_t}$ for \emph{strongly regular parameters} $t\in I_0=[\varepsilon_0,2\varepsilon_0]$: these are the parameters $t\in \bigcap\limits_{m=0}^{\infty} I_m$ where $I_0\supset I_1\supset\dots\supset I_m\supset\dots$ is a decreasing sequence of strongly regular intervals $I_m$. 

\begin{remark} It is interesting to notice that the strongly regular parameters of Palis-Yoccoz are not defined a priori, i.e., one has to perform the entire induction scheme before putting the hands on them. This is in contrast with the so-called \emph{Jakobson theorem} \cite{J}, a sort of $1$-dimensional version of Theorem \ref{t.PY}, where the strongly regular parameters are known since the beginning of the argument (because the location of the critical set is known in advance for $1$-dimensional endomorphisms). 
\end{remark}

Before starting the analysis of strongly regular parameters, one needs to ensure that such parameters exist, that is, one want to know whether there are parameters left from the parameter exclusion scheme of J. Palis and J.-C. Yoccoz. This issue is carefully treated in Section 9 of \cite{PY}, where the authors estimate the relative speed of strips associated to elements $(P,Q,n)\in\mathcal{R}(I)$ when the parameter $t\in I$ moves, and, by induction, they are able to control the measure of bad (not strongly regular) intervals: as it turns out, the measure of the set of bad intervals is $\leq \varepsilon_0^{1+\tau^2}$, so that the strongly regular parameters $t\in I_0=[\varepsilon_0,2\varepsilon_0]$ have almost full measure in $I_0$, i.e., $\geq\varepsilon_0(1 - \varepsilon_0^{\tau^2})$ (cf. Corollary 15 of \cite{PY}).

\begin{remark}\label{r.PYheteroclinic} In order to get some strongly regular parameter, one has to ensure that the initial interval $I_0$ is strongly regular (otherwise, one ends up by excluding $I_0$ in the first step of Palis-Yoccoz induction scheme, so that one has no parameters to play with in the next rounds of the induction). Here, J. Palis and J.-C. Yoccoz make use of the technical assumption that one is unfolding a heteroclinic tangency. The idea is that the formation of bicritical elements takes a long time in heteroclinic tangencies because the points in the critical region should pass near $p_s$ first, then near $p_u$ and only then they can return to the critical region again; of course, in the case of homoclinic tangencies, it may happen that bicritical elements pop up quickly and this is why one can't include homoclinic tangencies in the statement of Theorem \ref{t.PY}.
\end{remark}

From now on, let us fix $t\in \bigcap\limits_{m=0}^\infty I_m$ a strongly regular parameter, and let's study $\Lambda_g$ for $g=g_t$. Keeping this goal in mind, we introduce 
$\mathcal{R}=\mathcal{R}(t)=\bigcup\limits_{m=0}^{\infty}\mathcal{R}(I_m)$ the collection of all affine-like iterates of $g$ coming from the strongly regular intervals $I_m$. Using the class 
$\mathcal{R}$, we can define the class $\mathcal{R}_+^{\infty}$ of \emph{stable curves}, i.e., the class of curves $\omega$ coming from intersections $\omega=\bigcap\limits_{m=0}^\infty P_m$ of decreasing sequences $P_0\supset P_1\supset\dots$ of vertical strips serving as domain of affine-like iterates of $g$, that is, $(P_m,Q_m,n_m)\in\mathcal{R}$. Also, we put $\widetilde{\mathcal{R}}_+^{\infty}=\bigcup\limits_{\omega\in\mathcal{R}_+^{\infty}}\omega\subset M$ the set of points of $M$ in some stable curve.

These stable curves were introduced by analogy with uniformly hyperbolic horseshoes: indeed, the stable lamination of $K_g$ can be recovered from the transitions maps $g:R_a\cap g^{-1}(R_a)\to g(R_a)\cap R_{a'}$ by looking at the decreasing sequences of domains of simple compositions of these transitions maps. 

From the nice features of strong regular parameters, it is possible to prove that the class $\mathcal{R}_+^{\infty}$ is a $C^{1+Lip}$-lamination and one can use $g$ to induce a dynamical system $T^+:\mathcal{D}_+\subset \mathcal{R}_+^{\infty}\to\mathcal{R}_+^{\infty}$ isomorphic to a Bernoulli map with infinitely many branches (cf. Subsection 10.4 of \cite{PY}). Here, $\mathcal{D}_+$ is the set of stable curves not contained in infinitely many \emph{prime}\footnote{We say that an element $(P,Q,n)\in\mathcal{R}$ is \emph{prime} if it can't be obtained by simple composition of shorter elements $(P_0,Q_0,n_0), (P_1,Q_1,n_1)\in\mathcal{R}$ (shorter meaning $n_0,n_1<n=n_0+n_1$).} elements of 
$\mathcal{R}$. In particular, as it is shown in Subsections 10.5, 10.6, 10.7, 10.8, 10.9, 10.10 of \cite{PY}, $T^+$ is a non-uniformly hyperbolic dynamical system (in a very precise sense). Of course, by the symmetry between past and future (see Remark \ref{r.past-future}), one also has an analogous non-uniformly hyperbolic dynamical system $T^-:\mathcal{D}_-\subset \mathcal{R}_-^{\infty}\to\mathcal{R}_-^{\infty}$ on unstable curves, so that $\Lambda_g$ inherits a natural \emph{non-uniformly hyperbolic part} consisting of points whose $T^+$ and $T^-$ iterates never escape $\mathcal{R}_+^{\infty}$ and $\mathcal{R}_-^{\infty}$. 

Therefore, if we can show that the size of the sets of the points of $\Lambda_g$ escaping $\mathcal{R}_+^{\infty}$ or $\mathcal{R}_-^{\infty}$ is relatively small compared to the non-uniformly hyperbolic part of $\Lambda_g$, then we can say that $\Lambda_g$ is a \emph{non-uniformly hyperbolic horseshoe}. Here, J. Palis and J.-C. Yoccoz set up in Section 11 of \cite{PY} a series of estimates towards showing that the points of $\Lambda$ escaping $\mathcal{R}_+^{\infty}$ or $\mathcal{R}_-^{\infty}$ are \emph{exceptional}: for instance, they show Theorem \ref{t.PY-thm7} that the $2$-dimensional Lebesgue measure of $W^s(\Lambda_g)$ is zero because this property is true for the non-uniformly hyperbolic part of $\Lambda_g$ (by the usual hyperbolic theory) and the set of points of $\Lambda_g$ escaping $\mathcal{R}_+^{\infty}$ or $\mathcal{R}_-^{\infty}$ are rare in the sense that their $2$-dimensional Lebesgue measure contribute as an error term to the the non-uniformly hyperbolic part of 
$\Lambda_g$. 

At this point, our overview of Palis-Yoccoz induction scheme is complete. Closing this subsection and the first (survey) part of this text, we would like to make two comments. Firstly, as it is pointed out in page 14 of \cite{PY}, the philosophy that $\Lambda_g$ is constituted of a non-uniformly hyperbolic part and an exceptional set makes them expect that one could improve the information on the geometry of $W^s(\Lambda_g)$ or $\Lambda_g$. As it turns out, we will discuss in the second part of this text some recent results in this direction. Finally, condition \eqref{e.PYcondition} is not expected to be sharp by any means, but it seems that the strongly regular parameters are not sufficient to go beyond \eqref{e.PYcondition}, so that it is likely that one has to exclude further parameters in order to improve Theorem \ref{t.PY}.

\section{Part II -- a research announcement on non-uniformly hyperbolic horseshoes}

In what follows, we will consider the same setting of the article \cite{PY} of J. Palis and J.-C. Yoccoz, and we will discuss a recent improvement (obtained in collaboration with J. Palis and J.-C. Yoccoz) on Theorem \ref{t.PY-thm7} above. In particular, \emph{all} statements below concern the dynamics of $\Lambda_g$ where $g=g_t$ and $t$ is a strongly regular parameter in the sense of \cite{PY}: in other words, in the sequel, we do not have to exclude further parameters in order to get our (slightly improved) statements.

The main result of this part of the text is:

\begin{theorem}[C. Matheus, J. Palis and J.-C.Yoccoz \cite{MPY}]\label{t.MPY} The Hausdorff dimension of the stable and unstable sets $W^s(\Lambda_g)$ and $W^u(\Lambda_g)$ of the non-uniformly hyperbolic horseshoe $\Lambda_g$ is strictly smaller than $2$.
\end{theorem}

Logically, this result improves Theorem \ref{t.PY-thm7} because any subset of the compact $2$-dimensional manifold $M$ with Hausdorff dimension strictly smaller than $2$ has zero $2$-dimensional Lebesgue measure. 

The plan for the rest of this text is the following: in the next subsection we will prove Theorem \ref{t.MPY}, and in the final subsection we will make some comments on further results obtained in \cite{MPY}.

\subsection{Hausdorff dimension of the stable sets of non-uniformly hyperbolic horseshoes} Start by nicely decomposing the stable set $W^s(\Lambda_g)$. Using the notations of Subsection \ref{ss.PY-global-view}, we can write 
$$W^s(\Lambda_g):=\bigcup\limits_{n\geq 0} g^{-n}(W^s(\Lambda_g,\widehat{R})\cap R)$$
where $W^s(\Lambda_g,\widehat{R}):=\bigcap\limits_{n\geq 0} g^{-n}(\widehat{R})$. 

Since $g$ is a diffeomorphism, it follows from item (e) of Proposition \ref{p.Hausdorff-dimension} that 
$$\textrm{HD}(W^s(\Lambda_g))=\textrm{HD}(W^s(\Lambda_g,\widehat{R})\cap R)$$

Now, we separate $W^s(\Lambda_g,\widehat{R})\cap R$ into its good (non-uniformly hyperbolic) part and its exceptional part as follows:
$$W^s(\Lambda_g,\widehat{R})\cap R=:\bigcup\limits_{n\geq 0} (W^s(\Lambda_g,\widehat{R})\cap R\cap g^{-n}(\widetilde{\mathcal{R}}_+^{\infty})) \cup\mathcal{E}^+$$
In other words, the good part of $W^s(\Lambda_g,\widehat{R})\cap R$ consists of points passing by the nice set $\widetilde{\mathcal{R}}_+^{\infty}$ of stable curves and the exceptional set $\mathcal{E}^+$ is, \emph{by definition}, the complement of the good part. 

The set $\widetilde{\mathcal{R}}_+^{\infty}$ is the ``good'' (non-uniformly hyperbolic) part of the dynamics and hence it is not surprising that J. Palis and J.-C. Yoccoz showed in Section 10 of \cite{PY} that $\widetilde{\mathcal{R}}_+^{\infty}$ has Hausdorff dimension $1+d_s$ where $d_s$ is close to the stable dimension $d_s^0$ of the initial horseshoe $K$. 

Thus, the proof of Theorem \ref{t.MPY} is reduced\footnote{Here we're using item (b) of Proposition \ref{p.Hausdorff-dimension}.} to show that $\textrm{HD}(\mathcal{E}^+)<2$. 

Now, we follow the discussion of Section 11.7 of \cite{PY} to decompose $\mathcal{E}^+$ by looking at successive passages through \emph{parabolic cores} of strips. More concretely, given an element $(P,Q,n)\in\mathcal{R}$, we define the \emph{parabolic core} $c(P)$ of $P$ as
$$c(P)=\{p\in W^s(\Lambda_g,\widehat{R}): p\in P \textrm{ but }p\notin P' \textrm{ for all }P' \textrm{ \emph{child} of }P\}$$
Here, a child $P'$ of $P$ is a $(P',Q',n')\in\mathcal{R}$ such that $P'\subset P$ but there is no $P'\subset P''\subset P$ with $(P'',Q'',n'')\in\mathcal{R}$. The geometry of a parabolic core $c(P_k)$ of $(P_k,Q_k,n_k)\in\mathcal{R}$ is depicted below:

\begin{figure}[htb!]
\includegraphics[scale=0.53]{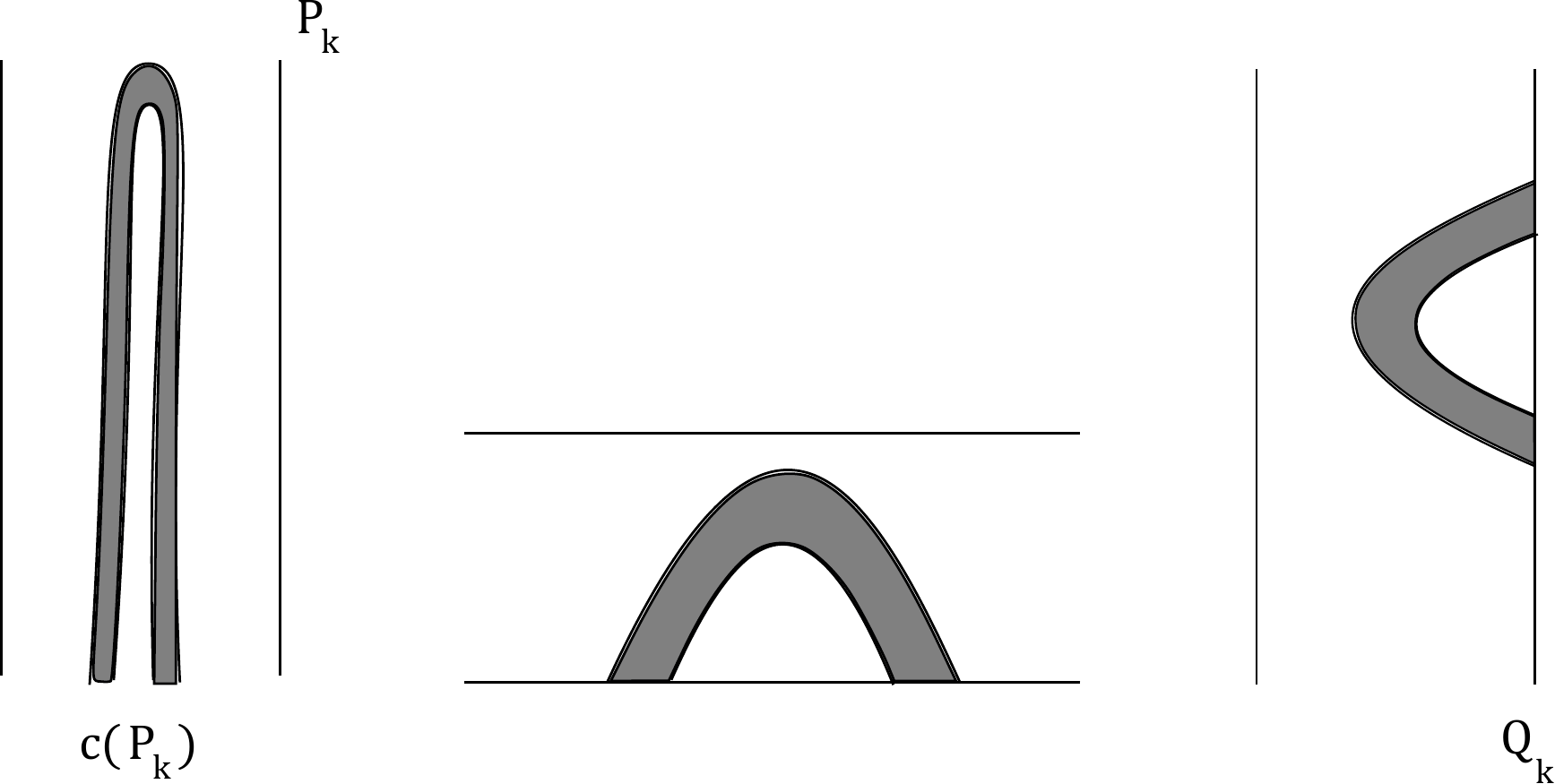}
\caption{The parabolic core $c(P_k)$ of $P_k$ belongs to the grey region inside $P_k$.}\label{f.16}
\end{figure}

By checking the definitions (of good part and exceptional set), it is not hard to convince oneself that $\mathcal{E}^+$ can be decomposed as
$$\mathcal{E}^+=\bigcup\limits_{(P_0,\dots,P_k) \textrm{ admissible}}\mathcal{E}^+(P_0,\dots,P_k)$$ 
where $(P_0,Q_0,n_0),\dots,(P_k,Q_k,n_k)\in\mathcal{R}$ and the sets $\mathcal{E}^+(P_0,\dots,P_k)$ are inductively defined as $\mathcal{E}^+(P_0):=c(P_0)$, $\mathcal{E}^+(P_0,P_1)=\{z\in\mathcal{E}^+(P_0): g^{n_0}(G(z))\in c(P_1)\}$, $\dots$ (cf. Equations (11.57) to (11.63) in \cite{PY}). Here, we say that $(P_0,\dots,P_k)$ is \emph{admissible} if $\mathcal{E}^+(P_0,\dots,P_k)\neq\emptyset$.

As the picture above indicates, the fact that the points of $\mathcal{E}^+(P_0,\dots,P_k)\neq\emptyset$ pass by successive parabolic cores imposes strong conditions over the elements $(P_i,Q_i,n_i)$: for instance, the parabolic core $c(P_i)$ of any $P_i$ is non-empty, the horizontal bands $Q_i$ are always critical and, because we're dealing $g=g_t$ where $t$ is a strongly $\beta$-regular parameter, the following estimate holds:
\begin{lemma}[Lemma 24 of \cite{PY}]\label{l.24}Suppose that $\mathcal{E}^+(P_0,\dots,P_{j+1})\neq\emptyset$, where $j\geq 1$. Then, 
$$\max(|P_{j+1}|,|Q_{j+1}|)\leq C|Q_j|^{\widetilde{\beta}}$$
where $\widetilde{\beta}=\beta(1-\eta)(1+\tau)^{-1}>1$.
\end{lemma}

This lemma is crucial for our purposes because it says that the exceptional set is confined 
into regions whose widths are decaying in a double exponential way 
to zero. Note that this is in sharp contrast with the case of the stable set of 
the initial horseshoe (which is confined into regions whose widths 
are going exponentially to zero): in other words, this lemma is a quantitative way of saying that the set $\mathcal{E}^+$ is exceptional when compared with the stable lamination of the horseshoe $K_g$.

In any event, the lemma above allow us to estimate the Hausdorff $d$-measure of $\mathcal{E}^+(P_0,\dots,P_k)$. By definition, we know that a certain $g$-iterate of 
$\mathcal{E}^+(P_0,\dots,P_k)$ is contained in the parabolic core $c(P_k)$. On the other hand, as it is shown in the figure below, we know that $c(P_k)$ is contained in a vertical strip of width $\varepsilon_k:=|Q_k|^{(1-\eta)/2}|P_k|$ and height $|Q_{k-1}|^{1/2}$ (see Proposition 62 of~\cite{PY}).

\begin{figure}[htb!]
\includegraphics[scale=0.53]{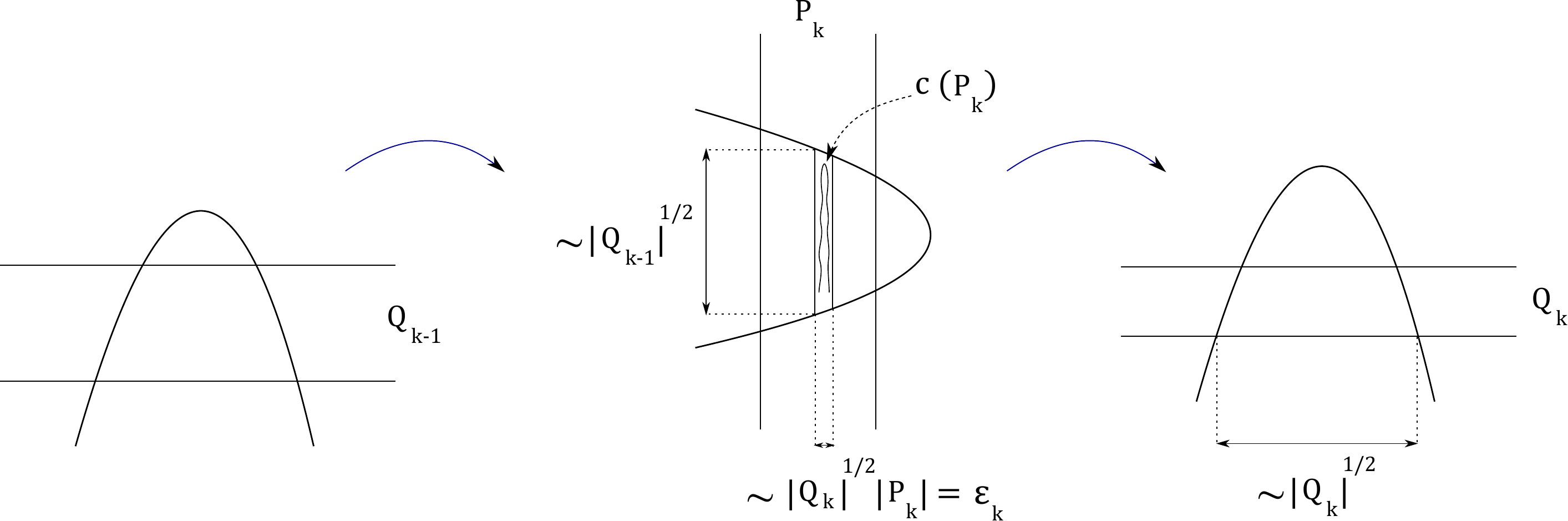}
\caption{Geometry of the parabolic core $c(P_k)$ of $P_k$.}\label{f.17}
\end{figure} 

We divide\footnote{Of course, this crude partition of $c(P_k)$ into $N_k$ squares of dimensions 
$\varepsilon_k\times\varepsilon_k$ aligned along a vertical strip is motivated by the fact that we do not want to keep track of the fine geometry of $\mathcal{E}^+(P_0,\dots,P_k)$ because it gets complicated very fast.} this vertical strip containing $c(P_k)$ into $N_k:=  \frac{|Q_{k-1}|^{1/2}}{|Q_k|^{(1-\eta)/2} |P_k|}$ squares of sides of lengths $\varepsilon_k$ and we analyze individually their evolution under the dynamics by an inductive procedure. More precisely, at the $i$-th step of our procedure, we have $N_{i+1}$ squares of dimensions $\varepsilon_{i+1}\times\varepsilon_{i+1}$ inside $Q_i$. We fix one of these squares and we note that $g^{-n_i}$ sends this square into a vertical strip of width $\varepsilon_i:=\varepsilon_{i+1}|P_i|$ and height $\varepsilon_{i+1}/|Q_i|:=\varepsilon_i/|P_i||Q_i|$ because $(P_i,Q_i,n_i)\in\mathcal{R}$ is affine-like. Again, we divide this vertical strip into $N_i:=1/|P_i| |Q_i|$ squares of sides of length $\varepsilon_i$ (similarly to Figure \ref{f.17}). Of course, during each step of this backward inductive procedure, we need to verify the compatibility condition $\varepsilon_{i+1}<|Q_i|$. In the present case, this compatibility condition is automatically satisfied in view of the estimate of Lemma~\ref{l.24}.

In particular, at the final step of this argument, we obtain a covering of $\mathcal{E}^+(P_0,\dots,P_k)$ by a collection of $N_0=N_k/\prod\limits_{i=0}^{k-1}|P_i| |Q_i|$ squares of sides of length $\varepsilon_0=\varepsilon_k\prod\limits_{i=0}^{k-1}|P_i|$. Thus, we have that 
\begin{eqnarray}\label{e.d-contribution-1}
\inf\limits_{\substack{\mathcal{O} \textrm{ cover of } \mathcal{E}^+(P_0,\dots,P_k) \\ \textrm{with } \textrm{diam}(\mathcal{O})<\varepsilon_0}} \, \, \sum\limits_{O_i\in\mathcal{O}} \textrm{diam}(O_i)^{d} &\leq& N_0\cdot\varepsilon_0^d = N_k\cdot\varepsilon_k^d\cdot\prod\limits_{i=0}^{k-1}|P_i|^{d-1}/|Q_i| \\  
&=& |Q_k|^{\frac{(1-\eta)(d-1)}{2}}\cdot |P_k|^{d-1}\cdot \frac{|P_{k-1}|^{d-1}}{ |Q_{k-1}|^{1/2}}\cdot\prod\limits_{i=0}^{k-2}\frac{|P_i|^{d-1}}{|Q_i|}\nonumber
\end{eqnarray}
Since $\widetilde{\beta}(1-\eta)(d-1)>1$ for any $d<2$ sufficiently close to $2$, we can use Lemma~\ref{l.24} to see that the right-hand side of \eqref{e.d-contribution-1} can be estimated by 
$$|Q_k|^{(1-\eta)(d-1)/2}\cdot |P_k|^{d-1}\cdot |Q_{k-1}|^{-1/2}\cdot |P_0|^{d-1}.$$
Let us take $d^->d_s^0+d_u^0-1$ a real number very close to $d_s^0+d_u^0-1$ and rewrite the previous expression as 
$$ |Q_k|^{\frac{(1-\eta)(d-1)}{2} - d^-}\cdot |Q_k|^{d^-} \cdot |P_k|^{d-1}\cdot |Q_{k-1}|^{-1/2}\cdot |P_0|^{d-1}.$$
Applying again Lemma~\ref{l.24}, we can bound this expression by 
$$ |Q_{k-1}|^{d^*}\cdot |Q_k|^{d^-} \cdot |P_0|^{d-1}$$
where $d^*=\frac{3(d-1)}{2}-d^- - 1/2$. However, the hypothesis (H4) of~\cite{PY} forces $0\leq d^-<1/5$, so that $d^*\geq 0$ for any $d\geq 22/15$.  It follows that 
\begin{equation}\label{e.d-contribution-2}
\inf\limits_{\substack{\mathcal{O} \textrm{ cover of } \mathcal{E}^+(P_0,\dots,P_k) \\ \textrm{with } \textrm{diam}(\mathcal{O})<\varepsilon_0}} \, \, \sum\limits_{O_i\in\mathcal{O}} \textrm{diam}(O_i)^{d}\leq |Q_{k-1}|^{d^*}\cdot |Q_k|^{d^-} \cdot |P_0|^{d-1}\leq |Q_k|^{d^-} \cdot |P_0|^{d-1}.
\end{equation}
Now we use two facts derived in the pages 204 and 205 of~\cite{PY}. Firstly, they show that the number of admissible sequences $(P_0,\dots,P_k)$ with fixed extremities $(P_0,Q_0,n_0)$ and $(P_k,Q_k,n_k)$ is $|Q_k|^{-\eta}$ (see Equation (11.77) of \cite{PY}). Secondly, the sum $\sum_{Q_k} |Q_k|^{d^- - \eta}$ is bounded because the results of the Subsection 11.5.9 of~\cite{PY} show that $\sum\limits_{Q \textrm{ critical }} |Q|^{d^- - \eta}$ converges\footnote{The fundamental fact that the critical locus is expected to have Hausdorff dimension $d_s+d_u-1$ is hidden in this estimate.} provided that $d^--\eta>d_s+d_u-1$. 

Putting these two facts together with~\eqref{e.d-contribution-2}, we see that the Hausdorff $d$-measure of $\mathcal{E}^+=\bigcup\limits_{(P_0,\dots,P_k) \textrm{ admissible}}\mathcal{E}^+(P_0,\dots,P_k)$ at scale $\varepsilon_0=\varepsilon_0(k)$ satisfies 
$$\inf\limits_{\substack{\mathcal{O} \textrm{ cover of } \mathcal{E}^+ \\ \textrm{with } \textrm{diam}(\mathcal{O})<\varepsilon_0}} \, \, \sum\limits_{O_i\in\mathcal{O}} \textrm{diam}(O_i)^{d}\leq \sum\limits_{P_0,\dots,P_k} |Q_k|^{d^-}\leq \sum\limits_{Q_k} |Q_k|^{d^- -C\eta}\leq C.$$  
Because $\varepsilon_0=\varepsilon_0(k)\to0$ as $k\to\infty$, this proves that $\textrm{HD}(\mathcal{E}^+)\leq d<2$. Hence, the proof of Theorem~\ref{t.MPY} is complete.

\subsection{Final comments on further results} The arguments of the previous subsection were based on \emph{soft analysis} of the geometry of the exceptional set. Every time the shape of the parabolic cores $c(P_i)$ was ready to get complicated, we divide it into squares and we analyzed the evolution of individual squares. In particular, every time we saw some parabolic geometry, we covered the ``tip of the parabola'' by a black box (square) and we forgot about the finer details of $\mathcal{E}^+$ in this region. Of course, it is not entirely surprising that this kind of soft estimate works to show $\textrm{HD}(W^s(\Lambda_g))<2$, but it is too crude if one wishes to compute the actual value of $\textrm{HD}(W^s(\Lambda_g))$. 

In particular, if one desires to prove that $\mathcal{E}^+$ is really exceptional so that $\textrm{HD}(W^s(\Lambda_g))$ is close to the expected dimension $1+d_s^0$, one has to somehow face the geometry of $\mathcal{E}^+$ and its successive passages through parabolic cores $c(P_i)$. 

\begin{figure}[htb!]
\includegraphics[scale=0.4]{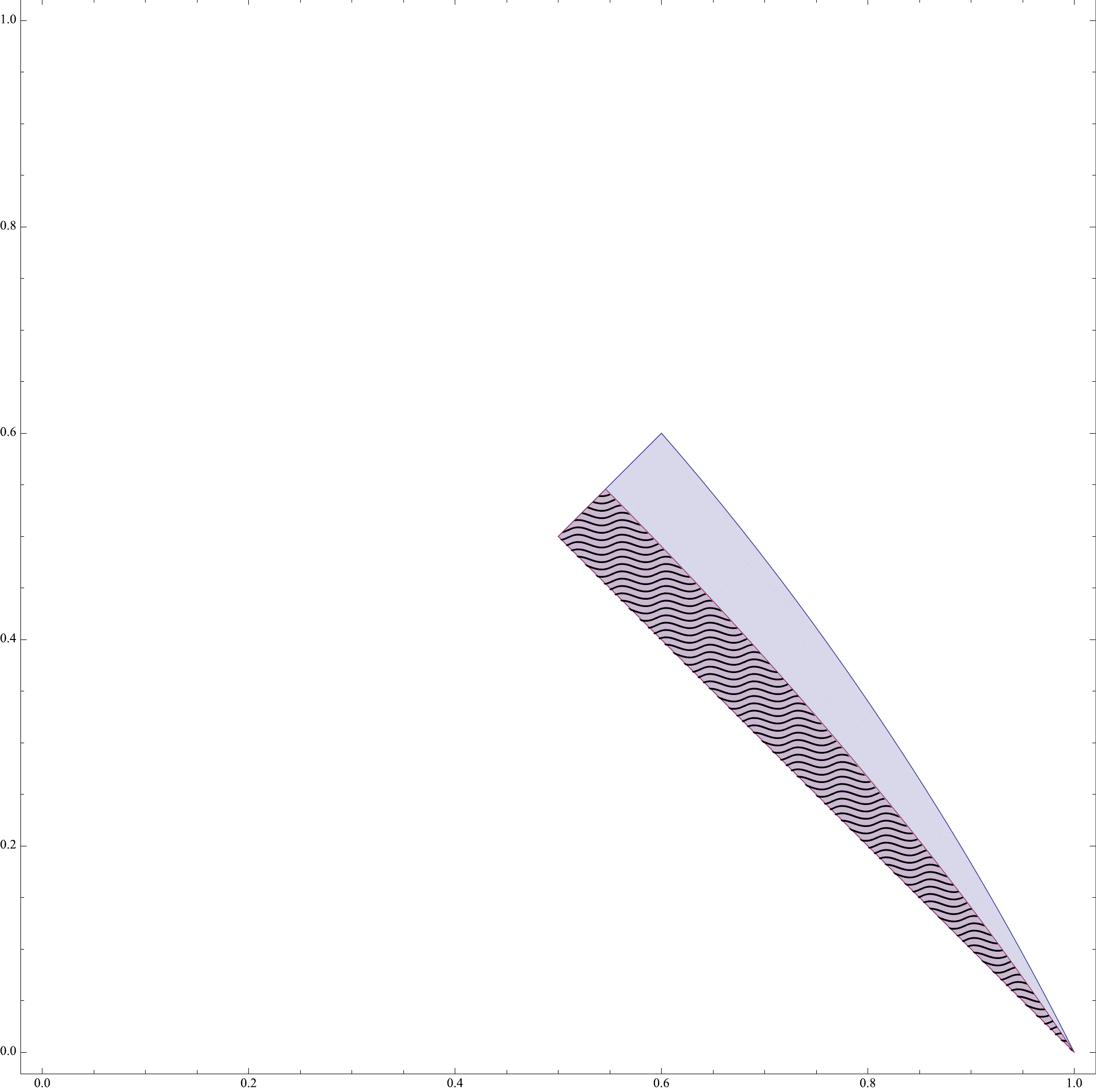}
\caption{$\mathcal{D}\cap\{d_s^0\leq d_u^0\}$ in wave texture sitting inside $\{(d_s^0,d_u^0): \eqref{e.PYcondition} \textrm{ holds}\}$.}\label{f.18}
\end{figure}

In the forthcoming article \cite{MPY}, we improve the soft strategy above without entering too much into the fine geometry of $\mathcal{E}^+$ by noticing that each $\mathcal{E}^+(P_0,\dots,P_k)$ is the image of the parabolic core $c(P_k)$ with a map of the form $g^{-n_0}\circ G\circ g^{-n_1}\circ\dots$ (obtained by alternating compositions of affine-like iterates of $g$ and the folding map $G$) whose derivative and Jacobian can be reasonably controlled. Using this control, we can study the Hausdorff $d$-measure of $\mathcal{E}^+$ at certain fixed scales $\varepsilon_0=\varepsilon_0(k)$ in terms of the geometry of $c(P_k)$ and the bounds on the derivative and the Jacobian of the map sending $c(P_k)$ into $\mathcal{E}^+$. By putting forward this estimate, we can show that $W^s(\Lambda_g)$ has the expected Hausdorff dimension (namely $1+d_s$) in a certain subregion $\mathcal{D}$ of values of stable and unstable dimensions $d_s^0$ and $d_u^0$ of the initial horseshoe. In Figure \ref{f.18} below we depicted in wave texture the region $\mathcal{D}\cap\{(d_s^0,d_u^0): d_s^0\leq d_u^0\}$ inside the larger region $\{(d_s^0,d_u^0): \textrm{Palis-Yoccoz condition }\eqref{e.PYcondition} \textrm{ holds}\}$. Actually, in this picture we drew only $\mathcal{D}\cap\{(d_s^0,d_u^0): d_s^0\leq d_u^0\}$ because the other half $\mathcal{D}\cap\{(d_s^0,d_u^0): d_s^0\geq d_u^0\}$ of $\mathcal{D}$ can be deduced by symmetry. 

The intersection $\mathcal{D}\cap\{(d_s^0,d_u^0):d_s^0=d_u^0\}$ of the region $\mathcal{D}$ with the diagonal (corresponding to the ``\emph{conservative case}''\footnote{The nomenclature ``conservative'' comes from the fact that the stable and unstable dimensions of any horseshoe of a area-preserving diffeomorphism coincide.}) can be explicitly computed (see \cite{MPY}):
$$\mathcal{D}\cap\{(d_s^0,d_u^0):d_s^0=d_u^0\}=\{(x,x): 1/2<x<0.545...\}$$
In particular, by putting this together with Figure \ref{f.18}, we see that $\mathcal{D}$ occupies slightly less than half of region given by Palis-Yoccoz condition \eqref{e.PYcondition}.

Finally, let us remark that we get the expected Hausdorff dimensions for $W^s(\Lambda_g)$ and $W^u(\Lambda_g)$ (in region $\mathcal{D}$), but the arguments can not be used to get the expected Hausdorff dimension for $\Lambda_g=W^s(\Lambda_g)\cap W^u(\Lambda_g)$. In fact, our constructions so far start from the future of $W^s(\Lambda_g)$ and the past of $W^u(\Lambda_g)$ where some geometric control is available, e.g., in the form of nice partitions, and then it tries to bring back the information, i.e., partitions, by analyzing the $g$-iterates used in our way back to the present time. Of course, this works if we deal separately with the past or the future. If we try to deal with both at the same time, we run into trouble because it is not obvious how the partitions coming from the future and the past intersect in the present time (due to the lack of transversality produced by $g$-iterates related to the folding map $G$). Evidently, the question of getting the expected Hausdorff dimension for 
$\Lambda_g$ is natural and interesting, and we hope to address this issue in \cite{MPY}. 

\bibliographystyle{amsplain}

\end{document}